\documentclass[11pt]{amsart}
\usepackage{graphicx}
\usepackage{amssymb}
\usepackage{amsmath} 
\usepackage{mathrsfs}
\usepackage{mathabx}
\usepackage{textcomp}
\usepackage{hyperref}
\usepackage{bbm}
\usepackage{pifont} 
\usepackage{color}
\usepackage{pgfplots}
\usepackage{cite}
\usepackage[margin=1in]{geometry}
\usepackage{subcaption}
\usepackage{slashbox}

\def\IN{\mathbb N}
\def\INo{{\mathbb N}_0}
\def\1I{\mathbbm{1}}
\def\tB{\mathtt{B}}
\def\th{\mathtt{h}}
\def\IZ{\mathbb Z}

\def\IP{\mathbb P}

\def\IE{\mathbb E}

\def\A{\mathcal A}
\def\B{\mathcal B}
\def\C{\mathcal C}

\def\F{\mathcal F}
\def\G{\mathcal G}

\def\L{\mathcal L}

\def\S{\mathcal S}
\def\T{\mathcal T}

\def\ms{\medskip\noindent}
\theoremstyle{plain}      
\theoremstyle{plain}      \newtheorem{corollary}{Corollary}
\theoremstyle{plain}      \newtheorem{proposition}{Proposition}
\theoremstyle{plain}      \newtheorem{claim}{Claim}
\theoremstyle{plain}      \newtheorem{theorem}{Theorem}
\theoremstyle{plain}      \newtheorem{example}{Example}
\theoremstyle{plain}

\bibstyle{plain}
\title[Dominant vertices and Attractors' Landscape]{Dominant vertices and attractors' landscape for Boolean networks}
\author{A. Espa\~na, W. Funez and E. Ugalde}
\address{Instituto de F\'isica, Universidad Aut\'onoma de San Luis 
Potos\'i, M\'exico.}
\email{andreae@ifisica.uaslp.mx, ugalde@ifisica.uaslp.mx, william.funez@if.uaslp.mx}
 
\date{}                                         
\begin{document}
\maketitle

\baselineskip=15pt

\centerline{\emph{Dedicated to the memory of Claudine Chaouiya}}

\begin{abstract}
In previous works, we introduced the notion of dominant vertices in the context of dynamical systems on networks. This is a set of nodes in the underlying network whose evolution determines the whole network's dynamics after a transient time.  In this paper, we focus on the case of Boolean networks. We define a reduced graph on the dominant vertices and an induced dynamics on this graph, which we prove is asymptotically equivalent to the original Boolean dynamics. Asymptotic conjugacy ensures that the systems, restricted to their respective attractors, are dynamically equivalent. For a significant class of networks, the induced dynamics is indeed a reduction of the original system. In these cases, the reduction, which is obtained from the structure of dominant vertices, supplies a more tractable system with the same structure of attractors as the original one. Furthermore, the structure of the induced system allows us to establish bounds on the number and period of the attractors, as well as on the reduction of the basin's sizes and transient lengths. We illustrate this reduction by considering a class of networks, which we call clover networks, whose dominant set is a singleton. To get insight into the structure of the basins of attraction of Boolean networks with a single dominant vertex, we complement this work with a numerical exploration of the behavior of a parametrized ensemble of systems of this kind. 

\end{abstract}

\section{Introduction}
\ms Boolean networks are discrete models, first proposed by Kauffman~\cite{Kauffman1969} in the context of genetic regulatory dynamics, to explain the origin and diversity of cellular types. Since then, the literature concerning this kind of model has proliferated, ranging from pure theoretical studies (starting with the classical contributions by Thomas, Derrida, Parisi, and coauthors~\cite{Thomas1973, DerridaPomeau1986, DerridaWeisbuch1986,BastollaParisi1997,BastollaParisiRelevant1998,BastollaParis1998}) to very concrete applications to model specific biological systems (see, for instance~\cite{EspinosaSotoEtal2004, SanchezThieffry1997} among many others). This paper aims to contribute to the study of the relation between some rough characteristics of the Boolean network and the structure of the landscape of attractors it generates. 

\ms
In~\cite{LunaUgalde2008}, we showed that in a dissipative dynamical system defined over directed networks, the whole dynamics is determined, after a transient time, by the projection of the orbits on a subset of nodes we called dominant vertices. As we argue below, the same effect takes place in a Boolean network, where, due to the finiteness of the local rules (which are typically many-to-one), the dynamics can in general be considered as dissipative. It is natural then to ask whether a reduced dynamical system can be defined on a network of the dominant vertices only, which would eventually contain all the relevant information of the system. In this paper, we define a reduced graph on the dominant vertices and an induced dynamics on it, and we prove that this induced system is asymptotically equivalent to the original Boolean dynamics. The asymptotic equivalence, which we introduce here, relates systems that are dynamically equivalent once restricted to their respective attractors. Under some conditions on the underlying network, the induced network is indeed a reduction of the original system. In this situation, the induced network is a more tractable version of the original Boolean network that preserves the structure of attractors. Furthermore, the nature of the induction allows us to establish bounds on the number and period of the attractors, as well as on the reduction of the basin's sizes and transient lengths, depending only on the structure of the dominant set. To illustrate all this notions and explore the relation between structure of dominant vertices and the landscape of attractors, we consider a class of Boolean networks with a single dominant vertex, which we call clover networks. We complement this work with a numerical exploration of the behavior of a parametrized ensemble of systems of this kind. 

\ms
\paragraph{\bf Related work}
In the context of continuous regulatory dynamics on networks, Mochizuki and coauthors~\cite{MochizukiEtal2013A, MochizukiEtal2013B} proved that any feedback vertex set is a control set, as well as a determining set. Their results are completely analogous to ours, which refers to piecewise contractions, and were published a few years before theirs~\cite{LunaUgalde2008}. Indeed, as we proved in Proposition~\ref{prop:CycleDominatSets} below, feedback vertex sets and dominant sets coincide for networks with everywhere nonzero input degree. Feedback vertex sets were used by Aracena~\cite{Aracena2008} to bound the number of steady states admitted by a Boolean network. This result was later refined by considering some particularities of the dynamics~\cite{AracenaRichardSalinas2017, Richard2018}. The existence of a subset of nodes determining the behavior of the whole Boolean network has been considered at least since 1988~\cite{FlyvbjergKjaer1988}, where Flyvbjerg and Kj{\ae}r rigorously study random Boolean networks with uniform input degree $=1$. 
Bastolla and Paris~\cite{BastollaParisiRelevant1998} extend that notion to random Boolean networks with uniform input degree $\geq 1$. Relevant nodes are those whose state changes affect the dynamics of the network. Their size and distribution depend on the topology of the network as well as the particularities of the regulatory rules, and are correlated to the size and distribution of the attractors. A notion similar to that of relevant nodes, this one referred to the controllability of Boolean networks, is that of driver nodes, introduced in~\cite{LiuEtal2011} in the context of continuous dynamics on complex networks, and further studied by Akutsu and coauthors in the case of Boolean networks~\cite{AkutsuEtal2007,HouEtal2016,HouEtal2019}. 

\ms The exponential growth of the configuration space in Boolean networks has motivated the search for more simplified versions of the original network while retaining all its important dynamical features. This has led to several model reduction proposals, from which we single out~\cite{NaldiEtal2011, Velizcuba2011,SaadatpourEtal2013}, where some transitional nodes are suppressed once their effect is taken into account by precomputing their possible output, and~\cite{ZanudoAlbert2013,ZanudoAlbert2015}, where sets of nodes defining stable motifs are assumed to attain a constant value and then replaced by their influence on the rest of the network. It is this kind of forward preprocessing that we use when defining the induced network in Section~\ref{sec:induced}. There are other unrelated reduction techniques, which consider other particularities of the dynamics, as for instance~\cite{ArgyrisEtal2023}, where the reduction is achieved by identification of backward equivalent nodes.

\ms
\paragraph{\bf The plan of the paper is as follows}
In Section~\ref{sec:preliminaries}, after reminding the definition of a Boolean network, we define the sets of dominant nodes, and propose a characterization of the attractors' landscape in terms of periods, basin sizes, and transient lengths. The main result of the section is the dominance-dynamics relation: if two trajectories coincide on a dominant set for a sufficiently long time, then they fully coincide thereafter. In Section~\ref{sec:induced} we define the induced automata network associated with a dominant set and establish the main results on eventual conjugacy, together with bounds on key dynamical indicators derived from the network topology.
Section~\ref{sec:CloverNetworks} is devoted to the class of Boolean networks with singleton dominant sets, which we call clover networks. We complement the section by presenting the results of a numerical exploration of ensembles of such networks, which serve as illustrations of the theoretical predictions. Finally, in Section~\ref{sec:Conclusions} we present some final comments and remarks.

\section{Preliminaries}\label{sec:preliminaries}

\noindent 
\paragraph{\bf The Boolean Network} Consider a directed graph with vertices in the finite set $V$, and arrows in the set $A\subset V\times V$, and for each $v\in V$, let $I(v):=\{u\in V: (u,v)\in A\}$ be the input set for $v$. We extend this notation to sets of nodes as follows: $I(U):=\bigcup_{u\in U}I(u)$ for each $U\subset V$. For each $s\in \INo$, define recursively $I^{s+1}(U)=I(I^s(U))$, with $I^0(U)\equiv U$. A Boolean network is a finite dynamical system defined on such a directed graph. To ensure the well-definiteness of the dynamics, we will always assume that the input set of every vertex $v\in V$ is non-empty. 

\ms Consider the set $\tB:=\{-1,1\}\subset\IZ$ and for each $x\in \tB^V$ and $U\subset V$, denote by $x_U\in \tB^U$ the restriction of the configuration $x$ to the vertices in $U$. In the case of a singleton, we will use the simpler notation $x_u$ instead of $x_{\{u\}}$. A Boolean network is a dynamical system on $\tB^V$, generated by the iteration of a function $F:\tB^V\to \tB^V$, with the following structure. For each $v\in V$, there exists a local rule $\phi_v:\tB^{I(v)}\to\tB$, such that $F(x)_v=\phi_v(x_{I(v)})$. We will refer to the graph $(V,A)$ used in the definition of $F$, as the underlying network. In order to simplify the notation, given $U\subset V$, we will use $\phi_U(x_{I(U)})$ to denote the configuration $\left(\phi_u\left(x_{I(u)}\right)\right)_{u\in U}\in \tB^U$. 

\ms Boolean networks belong to the larger class of automata networks (as the one considered in~\cite{Richard1986}, for instance). Automata networks are defined in a similar way as Boolean networks, with the difference that the states belong to an arbitrary finite set. Hence, an automata network with underlying graph $(V,A)$, taking values on the finite set $S$, is a dynamical system defined by the iteration of a function $F: S^V\to S^V$ such that, for each $v\in V$, there exists a function $\phi_v: S^{I(v)}\to S$, satisfying $F(x)_v=\phi_v(x_{I(v)})$.  

\ms 
\paragraph{\bf The attractors' landscape} The couple $(\tB^V, F)$ is a finite dynamical system whose behavior can be encompassed in a directed graph with vertices in $\tB^V$ and arrows $x\mapsto y$ for each couple $(x,y)\in\tB^V\times\tB^V$ such that $F(x)=y$. This directed graph $\T_{F}$ is the transition diagram of the dynamical system. It is easy to verify that each connected component $C$ of $\T_{F}$ consists of a single cycle decorated with a certain number of directed trees rooted at a point in the cycle. Each cycle codifies a periodic attractor of $(\tB^V, F)$, and the connected component containing this cycle codifies its entire basin of attraction. In this context, the attractors' landscape of $(\tB^V,F)$ is nothing but the structure of its transition diagram $\T_F$. We will characterize this structure by a list containing, for each connected component $C$, the length of the cycle $P(C)$ (i.e., the period of the attractor), the size of the component $|C|$ (i.e., the size of the basin of attraction), the mean length $\langle \tau\rangle_C$ and the maximal length $\tau_{\max}(C)$ of the attached trees.

\ms 
\paragraph{\bf Dominant nodes} Consider a network $ (V, A)$ such that $I(v)\neq \emptyset$ for each $v\in V$. We say that $U\subset V$ determines the set $U'\subset V$, if $I(U')\subset U$, i.e., if the input set of each node in $U'$ is completely contained in $U$. Given $U\subset V$, with $\partial U$ we denote the maximal set (in the sense of containment) determined by $U$, i.e., $\partial U=U' \Leftrightarrow (I(U')\subset U\, \wedge\, I(W)\subset U \Rightarrow W\subset U)$. The set $U\subset V$ is dominant if the chain $U_0:=U\subset U_1\subset\cdots\subset U_d:=V$, where $U_{\ell+1}=U_\ell\cup \partial U_\ell$ for each $0\leq \ell\leq d$. We call $U_0\subsetneq\cdots\subsetneq U_d$ the chain determined by $U$. Given a dominant set $U$, its chain is uniquely determined and we refer to its length $d=d(U)$, as the depth of the dominant set $U$. Alternatively, we can characterize the dominant sets is as follows. 

\ms
\begin{proposition}[Characterization of Dominant sets]~\label{prop:CycleDominatSets}
The set $U\subset V$ is dominant if and only if each cycle $v_0\mapsto v_1\mapsto\cdots \mapsto v_\ell\mapsto v_0$ in $(V,A)$ contains a vertex in $U$.  
\end{proposition}

\ms
The depth of the dominant set is nothing but the length of the longest path in $(V,A)$, starting at a vertex in $U$. The notion of a dominant set is closely related to the notion of feedback vertex set defined in~\cite{Aracena2008}. 

\begin{proof}\ 

\begin{itemize}
\item[$(\Rightarrow)$] Let $U_0\subsetneq\cdots\subsetneq U_d$ be the chain determined by $U$. Let us suppose that the cycle $v_0\mapsto \cdots \mapsto v_\ell\mapsto v_{\ell+1}=v_0$ does not contain vertices in $U$. If for a certain $i\geq 0$ we have $\{v_k:\ 0\leq k\leq \ell\}\cap U_i=\emptyset$, then $I(v_k)\setminus U_i\supset\{v_{k-1}\}\neq\emptyset$ for each $1\leq k\leq \ell+1$, and therefore $\{v_k:\ 0\leq k\leq \ell\}\cap \partial U_i=\emptyset$, which implies $\{v_k:\ 0\leq k\leq \ell\}\cap U_{i+1}=\emptyset$. It follows by induction that $\{v_k:\ 0\leq k\leq \ell\}\cap U_i=\emptyset$ for each $0\leq i\leq d$, which contradicts the hypothesis of $U$ being a dominant set. Hence, $v_0\mapsto \cdots \mapsto v_\ell\mapsto v_{\ell+1}=v_0$ necessarily contains a vertex in $U$. 

\item[$(\Leftarrow)$] If $U=V$, then $U_0:=U$ necessarily is a dominant set. Let us assume that $V\setminus U_0\neq \emptyset$. In this case there is a $v\in V\setminus U_0$ such that $I(v)\subset U_0$. Indeed, if $I(v)\setminus U_0\neq\emptyset$ for each $v\in V\setminus U_0$, then, by choosing $v_0\in V\setminus U_0$ and, for each $k\in\IN$, $v_{-k}\in I(v_{-k+1})\setminus U_0$, we obtain a path $v_{-|V|}\mapsto v_{-|V|+1}\mapsto\cdots \mapsto v_{0}$ such that $\{v_{-k}: 0\leq k\leq |V|\}\cap U_0=\emptyset$. This path, longer than the cardinality of $V$, must contain a cycle, contradicting the hypothesis. Therefore $\partial U_0\neq U_0$ and $U_1=U_0\cap \partial U_0\supsetneq U_0$. We recursively define $U_{i+1}:=U_i\cap \partial U_i$. Notice that, as long as $V\setminus U_i\neq \emptyset$, then, following the same reasoning as before, we conclude that there exists $v\in V\setminus U_i$ such that $I(v)\subset U_i$, and therefore $U_{i+1}=U_i\cap \partial U_i\supsetneq U_i$. Since $V$ is finite, this process must stop at some $d > 0$, and we conclude that $U$ is a dominant set with chain $U:=U_0\subsetneq\cdots\subsetneq U_d:=V$.

\end{itemize}

\end{proof}

\ms
\begin{example}  As an example, consider the network $(V,A)$ with $V=\{1,2,3,4,5\}$ and $A\subset V\times V$ as illustrate in  Figure~\ref{fig:Network2Dominats}. For this network, a dominant set of minimal size is $U=\{1,2\}$, determining the chain $U_0:=U\subset U_1:=\{1,2,3,4\}\subset U_2:=\{1,2,3,4,5\}:=V$. Hence, $d(U)=2$ in this case.

\begin{center}
\begin{figure}[h]
\begin{tikzpicture}
\draw[->] (2,3.7) -- (2,3.2);
\draw[->] (4,3.7) -- (4,3.2);
\draw[->] (4,2.7) -- (4,2.2);
\draw[->] (2.2,3.8) -- (3.7,3.2);
\draw[dashed,->] (1.5,3) to [out=180,in=180] (1.5,4);
\draw[dashed,->] (4.5,3) to [out=0,in=-0] (4.5,4);
\draw[->] (8.3,3) -- (9.7,3);
\draw[->] (7.8,2.8) .. controls (7,2.3) and (7,3.7) .. (7.8,3.2);
\draw[->] (10.2,3.2) .. controls (11,3.7) and (11,2.3) .. (10.2,2.8);

\filldraw [color=blue] (2,4) circle (3pt);
\filldraw [color=green] (2,3) circle (3pt);
\filldraw [color=blue] (4,4) circle (3pt);
\filldraw [color=green] (4,3) circle (3pt);
\filldraw [color=green] (4,2) circle (3pt);
\filldraw [color=blue] (8,3) circle (3pt);
\filldraw [color=blue] (10,3) circle (3pt);

\node at (1.7,4) {\text{1}};
\node at (1.7,3) {\text{3}};
\node at (4.3,4) {\text{2}};
\node at (4.3,3) {\text{4}};
\node at (4.3,2) {\text{5}};
\node at (8.2,2.7) {\text{1}};
\node at (9.7,2.7)  {\text{2}};

\end{tikzpicture}
\caption{A five-vertex network with a dominant set of two vertices colored in blue, and depth $d=2$. The reduced graph is shown aside.}~\label{fig:Network2Dominats}
\end{figure}
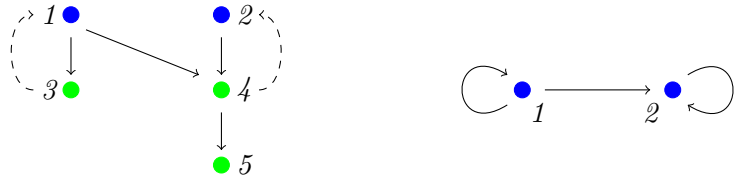   
\end{center}

\end{example}

\paragraph{\bf Remark on dominant sets}
The definition of dominant sets of vertices depends only on the topology of the underlying network and makes sense for Boolean networks with asynchronous updating. A given network admits multiple dominant sets: already, supersets of a dominant set are dominant as well. Even restricted to sets of minimal size, dominant sets are not unique. This can be easily verified by considering the directed cyclic graph, in which any one-vertex subset is dominant and of minimal size. Although dominant sets are easily characterized, identifying minimal dominant sets is, in general, computationally hard. In our contexts, finding a dominant set of minimal size is equivalent to finding a minimal feedback vertex set in a directed graph, a problem known to be NP-complete~\cite{Karp1972}.

\ms The following is an adaptation of a result established in~\cite{LunaUgalde2008}, to the present case. 

\begin{theorem}[Dominance and Dynamics]\label{theo:Dominance&Dynamics}
Consider the Boolean network $(\tB^V,F)$ with underlying graph  $(V,A)$. Let $U\subset V$ be a dominant set. If for $x\,,y\in \tB^V$, there exists $t_0\in\INo$ such that $F^{t}(x)_{U}=F^{t}(y)_U$ for $t_0\leq t\leq t_0+d(U)$, then $F^t(x)=F^t(y)$ for all $t\geq t_0+d(U)$. 
\end{theorem}

\ms According to this result, two orbits coinciding in the dominant set $U$ for an interval of time $t_0\leq t\, +\, t_0+d(U)$, necessarily collapse into the same orbit from time $t_0\,+\,d$. Hence, for two initial conditions $x,y\in \tB^V$ to belong to different basins of attraction, it is necessary that $F^t(x)_U\neq F^t(y)_U$ for all $t\in\INo$. Hence, a dominant set can be used to determine the attractors, their period, and their basins of attraction. The depth of the dominant set can also be used to bound the length of the transient behavior, as we will see below.

\ms 
\begin{proof}
The proof is quit direct. Since $x_U^{t_0}=y_U^{t_0}$,  then 
\[x^{t_0+1}_{\partial U}:=F(x)_{\partial U}=\phi_{\partial U}(x_U)=\phi_{\partial U}(y_U)=F(y)_{\partial U}=:y^{t_0+1}_{\partial U}.\]
Furthermore, since by hypothesis $x_U^{t_0+1}=y_U^{t_0+1}$, then we have $x_{U_1}^{t_0+1}=y_{U_1}^{t_0+1}$.  From here, we proceed by induction. Assuming that $x_{U_s}^{t_0+s}=y_{U_s}^{t_0+s}$, we obtain that 
\[x^{t_0+s+1}_{\partial U_s}:=F(x)_{\partial U_s}=\phi_{\partial U_s}(x_{U_s})=\phi_{\partial U_s}(y_U)=F(y)_{\partial U_s}=:y^{t_0+s+1}_{\partial U_s},\]
and since by hypothesis $x_U^{t_0+s+1}=y_U^{t_0+s+1}$, then we have $x_{U_{s+1}}^{t_0+s+1}=y_{U_{s+1}}^{t_0+s+1}$.  The result follows from the fact that $U_d=V$, which ensures that $x^{t_0+d}=y^{t_0+d}$, and therefore $x^t=y^t$ for all $t\geq t_0+d$.
\end{proof}

\bigskip
 
\section{The induced system} \label{sec:induced}

\ms
\paragraph{\bf The reduced graph} Consider the network $(V,A)$ with dominant set $U$ with depth $d=d(U)$, determining the chain $U_0:=U\subsetneq U_1\subsetneq\cdots\subsetneq U_d:=V$. The dominant set defines a reduced graph $(U,A_U)$, with $(u',u)\in A_U$ if and only if there exists a simple path in the original network, starting at $u'$ and ending at $u$.  

\ms For each $u\in U$, denote by $I_U(u)$ the input set of $u$ with respect to $(U,A_U)$. For $u\in U$ and $u'\in I_U(u)$, let $\ell(u',u)$ denote the set of all lengths of simple paths in $(V,A)$, starting at $u'$ and ending at $u$. With this, we define the recurrence length of $U$, given by $\ell=\ell(U):=\max\{\ell(u',u):\, u,u'\in U\}$. Notice that $\ell\leq d+1$.

\ms Let us recall that $I(W):=\bigcup_{w\in W}I(w)$ for each $W\subset V$ and for each $s\in \INo$,  
$I^{s+1}(W)=I(I^s(W))$, with $I^0(W)\equiv W$.  Using this notation, $(u,u')\in A_U$, if and only if there exists $0\leq s\leq d$ such that $u\in I^{s+1}(u')$. 

\ms 
\paragraph{\bf The induced automata network} A Boolean network $(\tB^V,F)$ with underlying graph $(V,A)$, induces an automata network on the reduced graph $(U,A_U)$. Let $\tB_\ell:=\tB^\ell$ and define $\Phi:\tB_\ell^U\to\tB$ as follows. For $y=\left(y^t_u\right)_{0\leq t < \ell;\, u\in U}\in \tB_\ell^U$, $\Phi(y):=\phi_{I(U)}(x^0)$, with $x^0\in \tB^{I(U)}$ given by
\begin{center}
\begin{tabular}{ll}
$x^0_{I(U)\cap U}= y^0_{I(U)\cap U}$, & $x^0_{I(U)\setminus U}= \phi_{I(U)\setminus U}(x^1)$ , 
\end{tabular}
\end{center}

\ms and for each $t\geq 1$, $x^t\in \tB^{I^{t+1}(U)}$ is defined recursively as
\begin{center}
\begin{tabular}{ll}
$x^t_{I^t(U)\cap U}=y^t_{I^t(U)\cap U}$, & $ x^t_{I^t(U)\setminus U}= \phi_{I^t(U)\setminus U}\left(x^{t+1}\right)$. 
\end{tabular}
\end{center}

\ms Since $I^{\ell-1}(U)\setminus U=\emptyset$, the recursion stops after $\ell$. With the collection $\Phi:\tB_\ell^{I_U}\to\tB$ constructed in this way, we define the automata network $(\tB_\ell^U, F)$, with $\F:\tB_\ell^U\to\tB^U$ given by
\begin{equation}\label{eq:InducedGLN}
\F(y)=\left(\Phi(y)\, y^0\,y^1\cdots y^{\ell-2}\right),
\end{equation}
for each $\left(y^0y^1\cdots y^{\ell-1}\right)\in\tB_d^U$. The configuration $y=\left(y^t_u\right)_{0\leq t< \ell;\, u\in U}\in \tB_\ell^U$ can be interpreted as the projection of a reversed orbit segment, $\left(F^{\ell-t}(x)\right)_{0\leq t<\ell}$, to the coordinates in the dominant set $U$. We will refer to this automata network as the automata network induced by $U$. 

\ms
\begin{example} Consider a Boolean network with underlying graph depicted in Figure~\ref{fig:Network2Dominats}, and defined by the collection of local functions $\{\phi_v:\tB^{I(v)}\to\tB\}_{v\in V}$. The dominant set $U=\{1,2\}$ of this underlying graph has recurrence length $\ell(U)=2$,  and the induced automata network $\F:\tB_2^{\{1,2\}}\to\tB_2^{\{1,2\}}$ is defined by the local functions:
\begin{center}
\begin{tikzpicture}
\draw[->] (1.3,3) -- (2.7,3);
\draw[->] (0.8,2.8) .. controls (0,2.3) and (0,3.7) .. (0.8,3.2);
\draw[->] (3.2,3.2) .. controls (4,3.7) and (4,2.3) .. (3.2,2.8);
\filldraw [color=blue] (1,3) circle (3pt);
\filldraw [color=blue] (3,3) circle (3pt);

\draw (1.2,2.7) node {\text{1}};
\node at (3,3.7) {$\Phi_1(y_1):=\left( \phi_1( \phi_3(y^1_1) )\,y_1^0\,\right)$};
\draw (2.7,2.7) node {\text{2}};
\node at (6,2.4) {$\Phi_2(y_1,y_2):=\left(\phi_2(\phi_4(y_1^1,y_2^1))\,y_2^0\,\right)$};

\end{tikzpicture}
\end{center}

\ms Notice how the definition of the local functions depends on the paths connecting dominant vertices to each other (see Figure~\ref{fig:Network2Dominats}).
\end{example}

\ms A configuration $x\in\tB^V$ is $F$-periodic of period $P\in\IN$, if $F^P(x)=x$. We denote by ${\rm Per}(F)$, the set of all $F$-periodic configurations, and with ${\rm Per}_P(F)$ we denote the subset of ${\rm Per}(F)$ containing all the configurations of period $p$. The minimal period for $x\in{\rm Per}(F)$ is the smallest $p\in\IN$ such that $x\in {\rm Per}_P(F)$. Similar notions apply for $\F$ acting on $\tB_d^U$. As a consequence of the previous result, we have the following.
 
\ms
\begin{theorem}[Eventual conjugacy]~\label{theo:Eventual-conjugacy}
Consider the Boolean network $(\tB^V,F)$ with underlying graph $(V,A)$. Let $U\subset V$ be a dominant set with recurrence length $\ell=\ell(U)$, and let $(\tB_\ell^U,\F)$ be the automata network induced by $U$. The transformation $\th:\tB^V\to \tB_\ell^U$, given by 
\begin{equation}\label{eq:Eventual-conjugacy}
\th(x)=\left(F^{\ell-1}(x)_U\,F^{\ell-2}(x)_U\cdots F(x)_U\,x_U\right),
\end{equation}
is such that 
\begin{itemize}
\item[a)] $\F(\th(x))=\th(F(x))$ for all $x\in\tB^V$, and
\item[b)] $\th(x)=\th(x')$ implies $x=x'$, for all $x,x'\in {\rm Per}(F)$.
\end{itemize}
\end{theorem}

\ms 
\begin{proof}
For each $y=(y^0y^1\cdots y^{\ell-1})\in \tB_d^U$, let us denote $(\Phi_u(y_{I_U(u)}))_{u\in U}\in \tB^U$ by $\Phi(y)$. With this, $\F\left(y^0y^1\cdots y^{\ell-1}\right)=\left(\Phi(y) y^0y^1\cdots y^{\ell-2}\right)$ and therefore,
\[
\F(\th(x))=\left(\Phi(\th(x))\,F^{\ell-1}(x)_U\,F^{\ell-2}(x)_U\cdots F(x)_U\right).
\]
Claim a) follows from the fact that $\Phi(\th(x))=F^{\ell}(x)_U$, which we now establish. 

\ms Let us recursively define the configurations $x^0\in \tB^I(U),x^1\in\tB^{I^2(U)},\ldots,x^{\ell-1}\in\tB^{I^\ell(U)}$, as follows: 

\begin{center}

\begin{tabular}{ll}
$x^0_{I(U)\cap U}=F^{\ell-1}(x)_{I(U)\cap U}$, & $x^0_{I(U)\setminus U}=\phi_{I(U)\setminus U}(x^1)$,\\
$x^1_{I^2(U)\cap U}=F^{\ell-2}(x)_{I^2(U)\cap U}$, & $x^1_{I^2(U)\setminus U}=\phi_{I^{2}(U)\setminus U}\left(x^{2}\right)$, \\
\vdots & \vdots   \\
$x^t_{I^{t+1}(U)\cap U}=F^{\ell-(t+1)}(x)_{I^{t+1}(U)\cap U}$, & $x^t_{I^{t+1}(U)\setminus U}=\phi_{I^{t+1}(U)\setminus U}\left(x^{t+1}\right)$, \\
\vdots & \vdots   \\
 $x^{\ell-2}_{I^{\ell-1}(U)\cap U}=F(x)_{I^{\ell-1}(U)\cap U}$, & $x^{\ell-2}_{I^{\ell-1}(U)\setminus U}=\phi_{I^{\ell-1}(U)\setminus U}\left(x^{\ell-1}\right)$,\\
  $x^{\ell-1}_{I^{\ell}(U)}=x_{I^{\ell}(U)}$. &  
\end{tabular}
\end{center}
Notice that, according to the definition of the induced automata network, $x^0$ defined in this way is such that $\phi_{I(U)}\left(x^0\right)=\Phi(y)$ with $y=\th(x)=\left(F^{\ell-1}(x)_U\,F^{\ell-2}(x)_U\,\cdots F(x)_U\,x_U\right)\in \tB_d$.

\ms 
We now prove that $\Phi(y)=\phi_{I(U)}\left(x^0\right)=F^{\ell}(x)_U$. For this, first note that $x^{\ell-1}_{I^{\ell}(U)}=x_{I^{\ell}(U)}\equiv F^0(x)_{I^{\ell}(U)}$. Now, assuming $x^{\ell-t}=F^{t-1}(x)_{I^{\ell-(t-1)}(U)}$, it follows that
\[x^{\ell-(t+1)}_{I^{\ell-t}(U)\setminus U}=\phi_{I^{\ell-t}(U)\setminus U}\left(F^{t-1}(x)_{I^{\ell-(t-1)}(U)}\right)
                                           =F^t\left(x\right)_{I^{\ell-t}(U)\setminus U}, \]
and since by the above definition $x^{\ell-(t+1)}_{I^{\ell-t}(U)\cap U}= F^t\left(x\right)_{I^{\ell-t}(U)\cap U}$, then $x^{\ell-(t+1)}=F^t(x)_{I^{\ell-t}(U)}$. In this way, we have inductively obtained $x^{\ell-t}=F^{t-1}(x)_{I^{\ell-(t-1)}(U)}$ for each $1\leq t\leq  \ell$, and in particular $x^0=F^{\ell}(x)_{U}$, and hence
\[\Phi(y)=\phi_{I(U)}\left(x^0\right)=F^{\ell}(x)_U,\]
and Claim a) follows.

\ms For claim b), consider $x, x'\in  {\rm Per}(F)$ such that $\th(x)=\th(x')$, i.e., $F^t(x)_U=F^t(x')_U$ for $0\leq t\leq d$. Theorem~\ref{theo:Dominance&Dynamics} implies that $F^t(x)=F^t(x')$ for all $t\geq d$, and since $F^p(x)=x$, $F^q(x')=x'$ for some $p,q\in\IN$, then necessarily
\[x=F^{\kappa\,pq}(x)=F^{\kappa\,pq}(x')=x',\]
for all $\kappa\in\IN$ such that $\kappa\,pq\geq d$, and claim b) follows.

\end{proof}

\ms
\paragraph{\bf Remarks to Theorem~\ref{theo:Eventual-conjugacy}}\

\begin{enumerate} 
\item Since the dominant set is not unique, the induced network could in principle depend on the choice of the dominant set. Natural choices of dominant sets are dominant sets of minimal cardinality, but even in that case, two dominant sets with the same minimal cardinality could have different recurrence lengths and non-equivalent reduced networks. 
 
\item According to this theorem, the induced automata network and the original Boolean network are semi-conjugate as dynamical systems. Furthermore, after a transient stage, the two systems are equivalent and, therefore have the same number and type of attractors. It is this equivalence in the long run, that we refer to as eventual conjugacy. A direct consequence of Theorem~\ref{theo:Eventual-conjugacy} is the fact that the transition diagrams $\T_F$  and $\T_\F$, defined by the Boolean network $(\tB^V, F)$ and the induced automata network $(\tB_d^U,\F)$ respectively, have the same number of connected components and that $\th$ defines a bijection between the cycles at the core of the connected components. Hence, the difference between $\T_F$ and $\T_\F$ lies in the structure of the trees rooted at the cycles at the core of each connected component. These trees codify the transient dynamics of each system, and the fact that the correspondence between $\T_F$ and $\T_\F$ is established through a semi-conjugacy, which gives place to a homomorphism between directed graphs, imposes some restrictions as we will see below.

\item The eventual conjugacy is an equivalence relation between dynamical systems on finite sets. Reflexivity and transitivity are obvious. To prove symmetry, consider the following: 

\begin{claim}
Let $(\A,f)$ and $(\B,g)$ be finite dynamical systems and $\th:\A\to\B$ is an eventual equivalence. Then there exist and eventual equivalence $\th':\B\to\A$.  
\end{claim} 

\begin{proof}
Let us start by partitioning $\B=\bigsqcup_{n=1}^N\B_n$, as a disjoint union of basins. Each $\B_n$ consisting of initial conditions whose orbits conver to the same periodic attractor. 
For each $1\leq n\leq N$, let $\L_n:=\{b\in \B_n:\, g^{-1}\{b\}=\emptyset\}$, i.e.,  $\L_n$ are the leave of the connected component of the transition diagram $\T_\B$, inside $\B_n$.  If $\L_n=\emptyset$, then $(\B_n,g)$ is a periodic orbit and $\th$ is invertible on $\B_n$. In this case $\th'|_{\B_n}=\th^{-1}$. If on the contrary $\L_n\neq\emptyset$, take $t_0:=\min\{t\geq 0:\, \th^{-1}(g^{t_0}(b))\neq\emptyset\,\forall b\in\L_n\}$, and for each $b\in\L_n$ choose $a:=a(b)\in \th^{-1}(g^{t_0}(b))$. With this, define $\th'|_{\L_n}$ by $\th'(\G^t(b)):=\F^t(a(b))$ for each $b\in\L_n$. Since $\B_n=\bigcup_{t\geq 0}\G(\L_n)$, then $\th'|{\B_n}$ is well defined, and it is a eventual conjugacy by construction.

\end{proof}

\end{enumerate}
\begin{example} As an example, consider the Boolean network with  the  graph  depicted in Figure~\ref{fig:Network2Dominats} as underlying graph, and defined by the local functions

\begin{align}~\label{eq:LocalBoolean}
\phi_1(x_3)&=-x_3, \, \phi_2(x_4)=-x_4, \, \phi_3(x_1)=x_1,\\
\phi_4(x_4)&=x_1\times x_2, \,  \phi_5(x_4)=x_4. \nonumber 
\end{align}

\ms In this case $U=\{1,2\}$ is the dominant set, $\ell(U)=1$,  and the  induced  automata network $\F:\tB_2^{\{1,2\}}\to\tB_2^{\{1,2\}}$ has local functions
\begin{align}~\label{eq:LocalInduced}
\Phi_1(y_1)&:=\left( \phi_1( \phi_3(y^1_1) )\,y_1^0\,\right)=\left(\left(-y_1^1\right)\,y_1^0\right), \\
\Phi_2(y_1,y_2)&:=\left(\phi_2(\phi_4(y_1^1,y_2^1))\,y_2^0\,\right)=\left(\left(-y_1^1\times y_2^1\right)\,y_2^0\right).
\nonumber
\end{align}

\ms In Figure~\ref{fig:TransitionDiagrams} we show the transition diagrams $\T_F$, corresponding to the Boolean network $\left(\tB^{\{1,\ldots,5\}}, F\right)$ defined by the local functions in Equations~\eqref{eq:LocalBoolean}, and the transition diagram $\T_\F$, corresponding induced automata network $\left(\tB_2^2,\F\right)$ defined in Equations~\eqref{eq:LocalInduced}. The nodes of $\T_F$ are enumerated following the lexicographic order in $\tB^{´\{1,\ldots 5\}}$, while the nodes in $\T_F$ are couples in $\tB_2^{\{1,2\}}$, enumerating the elements in $\tB_2$ lexicographically.

\begin{center}
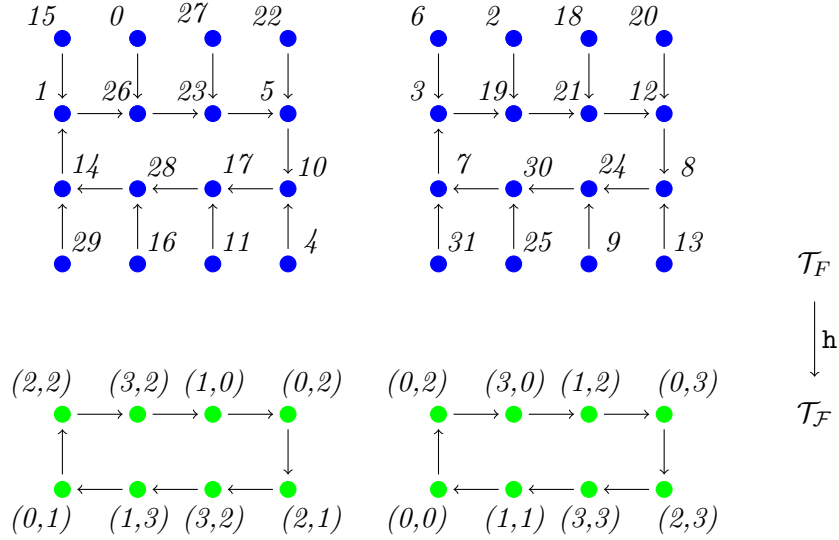
\begin{figure}[h]
\begin{tikzpicture}
\filldraw [color=blue] (1,8) circle (3pt);
\filldraw [color=blue] (2,8) circle (3pt);
\filldraw [color=blue] (3,8) circle (3pt);
\filldraw [color=blue] (4,8) circle (3pt);
\filldraw [color=blue] (1,7) circle (3pt);
\filldraw [color=blue] (2,7) circle (3pt);
\filldraw [color=blue] (3,7) circle (3pt);
\filldraw [color=blue] (4,7) circle (3pt);
\filldraw [color=blue] (1,6) circle (3pt);
\filldraw [color=blue] (2,6) circle (3pt);
\filldraw [color=blue] (3,6) circle (3pt);
\filldraw [color=blue] (4,6) circle (3pt);
\filldraw [color=blue] (1,5) circle (3pt);
\filldraw [color=blue] (2,5) circle (3pt);
\filldraw [color=blue] (3,5) circle (3pt);
\filldraw [color=blue] (4,5) circle (3pt);
\node at (0.7,8.3) {\text{15}};
\node at (1.7,8.3) {\text{0}};
\node at (2.7,8.3) {\text{27}};
\node at (3.7,8.3) {\text{22}};
\node at (0.7,7.3) {\text{1}};
\node at (1.7,7.3) {\text{26}};
\node at (2.7,7.3) {\text{23}};
\node at (3.7,7.3) {\text{5}};
\node at (1.3,6.3) {\text{14}};
\node at (2.3,6.3) {\text{28}};
\node at (3.3,6.3) {\text{17}};
\node at (4.3,6.3) {\text{10}};
\node at (1.3,5.3) {\text{29}};
\node at (2.3,5.3) {\text{16}};
\node at (3.3,5.3) {\text{11}};
\node at (4.3,5.3) {\text{4}};
\draw[->] (1,7.8) -- (1,7.2);
\draw[->] (2,7.8) -- (2,7.2);
\draw[->] (3,7.8) -- (3,7.2);
\draw[->] (4,7.8) -- (4,7.2);
\draw[->] (1,5.2) -- (1,5.8);
\draw[->] (2,5.2) -- (2,5.8);
\draw[->] (3,5.2) -- (3,5.8);
\draw[->] (4,5.2) -- (4,5.8);
\draw[->] (1,6.2) -- (1,6.8);
\draw[->] (4,6.8) -- (4,6.2);
\draw[->] (1.2,7) -- (1.8,7);
\draw[->] (2.2,7) -- (2.8,7);
\draw[->] (3.2,7) -- (3.8,7);
\draw[->] (3.8,6) -- (3.2,6);
\draw[->] (2.8,6) -- (2.2,6);
\draw[->] (1.8,6) -- (1.2,6);
\filldraw [color=blue] (6,8) circle (3pt);
\filldraw [color=blue] (7,8) circle (3pt);
\filldraw [color=blue] (8,8) circle (3pt);
\filldraw [color=blue] (9,8) circle (3pt);
\filldraw [color=blue] (6,7) circle (3pt);
\filldraw [color=blue] (7,7) circle (3pt);
\filldraw [color=blue] (8,7) circle (3pt);
\filldraw [color=blue] (9,7) circle (3pt);
\filldraw [color=blue] (6,6) circle (3pt);
\filldraw [color=blue] (7,6) circle (3pt);
\filldraw [color=blue] (8,6) circle (3pt);
\filldraw [color=blue] (9,6) circle (3pt);
\filldraw [color=blue] (6,5) circle (3pt);
\filldraw [color=blue] (7,5) circle (3pt);
\filldraw [color=blue] (8,5) circle (3pt);
\filldraw [color=blue] (9,5) circle (3pt);
\node at (5.7,8.3) {\text{6}};
\node at (6.7,8.3) {\text{2}};
\node at (7.7,8.3) {\text{18}};
\node at (8.7,8.3) {\text{20}};
\node at (5.7,7.3) {\text{3}};
\node at (6.7,7.3) {\text{19}};
\node at (7.7,7.3) {\text{21}};
\node at (8.7,7.3) {\text{12}};
\node at (6.3,6.3) {\text{7}};
\node at (7.3,6.3) {\text{30}};
\node at (8.3,6.3) {\text{24}};
\node at (9.3,6.3) {\text{8}};
\node at (6.3,5.3) {\text{31}};
\node at (7.3,5.3) {\text{25}};
\node at (8.3,5.3) {\text{9}};
\node at (9.3,5.3) {\text{13}};
\draw[->] (6,7.8) -- (6,7.2);
\draw[->] (7,7.8) -- (7,7.2);
\draw[->] (8,7.8) -- (8,7.2);
\draw[->] (9,7.8) -- (9,7.2);
\draw[->] (6,5.2) -- (6,5.8);
\draw[->] (7,5.2) -- (7,5.8);
\draw[->] (8,5.2) -- (8,5.8);
\draw[->] (9,5.2) -- (9,5.8);
\draw[->] (6,6.2) -- (6,6.8);
\draw[->] (9,6.8) -- (9,6.2);
\draw[->] (6.2,7) -- (6.8,7);
\draw[->] (7.2,7) -- (7.8,7);
\draw[->] (8.2,7) -- (8.8,7);
\draw[->] (8.8,6) -- (8.2,6);
\draw[->] (7.8,6) -- (7.2,6);
\draw[->] (6.8,6) -- (6.2,6);
\node at (11,5) {$\T_F$};
\draw[->] (11,4.5) -- (11,3.5);
\node at (11.2,4) {$\th$};
\node at (11,3) {$\T_\F$};
\filldraw [color=green] (1,3) circle (3pt);
\filldraw [color=green] (2,3) circle (3pt);
\filldraw [color=green] (3,3) circle (3pt);
\filldraw [color=green] (4,3) circle (3pt);
\filldraw [color=green] (1,2) circle (3pt);
\filldraw [color=green] (2,2) circle (3pt);
\filldraw [color=green] (3,2) circle (3pt);
\filldraw [color=green] (4,2) circle (3pt);
\filldraw [color=green] (6,3) circle (3pt);
\filldraw [color=green] (7,3) circle (3pt);
\filldraw [color=green] (8,3) circle (3pt);
\filldraw [color=green] (9,3) circle (3pt);
\filldraw [color=green] (6,2) circle (3pt);
\filldraw [color=green] (7,2) circle (3pt);
\filldraw [color=green] (8,2) circle (3pt);
\filldraw [color=green] (9,2) circle (3pt);
\node at (0.7,3.4) {\text{(2,2)}};
\node at (2,3.4) {\text{(3,2)}};
\node at (3,3.4) {\text{(1,0)}};
\node at (4.3,3.4) {\text{(0,2)}};
\node at (4.3,1.6) {\text{(2,1)}};
\node at (3,1.6) {\text{(3,2)}};
\node at (2,1.6) {\text{(1,3)}};
\node at (0.7,1.6) {\text{(0,1)}};
\node at (5.7,3.4) {\text{(0,2)}};
\node at (7,3.4) {\text{(3,0)}};
\node at (8,3.4) {\text{(1,2)}};
\node at (9.3,3.4) {\text{(0,3)}};
\node at (9.3,1.6) {\text{(2,3)}};
\node at (8,1.6) {\text{(3,3)}};
\node at (7,1.6) {\text{(1,1)}};
\node at (5.7,1.6) {\text{(0,0)}};
\draw[->] (1,2.2) -- (1,2.8);
\draw[->] (4,2.8) -- (4,2.2);
\draw[->] (1.2,3) -- (1.8,3);
\draw[->] (2.2,3) -- (2.8,3);
\draw[->] (3.2,3) -- (3.8,3);
\draw[->] (3.8,2) -- (3.2,2);
\draw[->] (2.8,2) -- (2.2,2);
\draw[->] (1.8,2) -- (1.2,2);
\draw[->] (6,2.2) -- (6,2.8);
\draw[->] (9,2.8) -- (9,2.2);
\draw[->] (6.2,3) -- (6.8,3);
\draw[->] (7.2,3) -- (7.8,3);
\draw[->] (8.2,3) -- (8.8,3);
\draw[->] (8.8,2) -- (8.2,2);
\draw[->] (7.8,2) -- (7.2,2);
\draw[->] (6.8,2) -- (6.2,2);
\end{tikzpicture}
\caption{Above, the transition diagram $\T_F$, corresponding to the Boolean network $(\tB^{\{1,\ldots,5\}},F)$. Below, the transition diagram $\T_\F$, defined by the induced logic network $(\tB_2^{\{1,2\}},\F)$. The eventual conjugacy $\th$, is such that $\th(1)=(2,2)$ and $\th(3)=(2,0)$, which establishes the equivalence between the two periodic attractors of $(\tB^{\{1,\ldots,5\}},F)$ with those of $(\tB_2^{\{1,2\}},\F)$.}~\label{fig:TransitionDiagrams}
\end{figure}
\end{center}
\end{example}

\ms In the previous example, the transition diagram of the induced automata network preserves the two cycles of the original Boolean network but has lost the ramifications that codify the transient dynamics. Below, we show another example illustrating how the eventual conjugacy between a Boolean network and its induced automata network, while affecting the ramifications rooted in the cycles of the transition diagram, partially preserves the transient regimes.

\ms
\begin{example}
Consider the Boolean network $\left(\tB^{\{1,\ldots,5\}},F\right)$ depicted in Figure~\ref{fig:Penacho}, defined by the local functions
\begin{align}~\label{eq:LocalBoolean2}
\phi_1(x_2,x_3)&=x_2\times x_3, \, \phi_2(x_1)=x_1, \, \phi_3(x_1)=-x_1,\\
\phi_4(x_2,x_3)&=-x_2\times x_3, \,  \phi_5(x_4)=-x_4. \nonumber 
\end{align}
For this network,  the dominant set is $U=\{1\}$, with depth $d=3$, and recurrence length $\ell=2$. Hence, the induced automata network is a transformation $\F:\tB_2\to\tB_2$, defined by the local function
\begin{equation}~\label{eq:LocalInduced2}
\Phi\left(y^0\,y^1\right)=\left(-y^1\times y^1\right)=-1.
\end{equation}
In the same figure, we show the transition diagrams of both the original Boolean network and the induced automata network. The eventual conjugacy $\th:\tB^{\{1,\ldots,5\}}\to\tB_2$ is codified by the correspondence of node colors in both diagrams. 

\begin{center}
\begin{figure}
\begin{tikzpicture}[scale=1.3]
\filldraw [color=blue] (4,9) circle (3pt);
\filldraw [color=blue] (3,10) circle (3pt);
\filldraw [color=blue] (5,10) circle (3pt);
\filldraw [color=green] (4,11) circle (3pt);
\filldraw [color=blue] (4,8) circle (3pt);
\draw (3.7,11) node {\text{1}};
\draw (2.7,10) node {\text{2}};
\draw (5.2,10) node {\text{3}};
\draw (3.7,9) node {\text{4}};
\draw (3.7,8) node {\text{5}};
\filldraw [color=green] (9,10) circle (3pt);
\draw (8.7,10) node {\text{1}};
\draw[->] (3.8,10.8) -- (3.2,10.2);
\draw[->] (4.2,10.8) -- (4.8,10.2);
\draw[->] (3.2,9.8) -- (3.8,9.2);
\draw[->] (3.8,10.8) -- (3.2,10.2);
\draw[->] (4.8,9.8) -- (4.2,9.2);
\draw[->] (4,8.8) -- (4,8.2);
\draw[dashed,->] (3,10.2) to [out=90,in=224] (3,11.5) to [out=45,in=135] (3.8,11.2);
\draw[dashed,->] (5,10.2) to [out=90,in=315] (5,11.5) to [out=115,in=45] (4.2,11.2);
\draw[->] (9.2,10) .. controls (9.7,10.1) and (9.3,10.7) .. (9,10.2);
\filldraw [color=blue] (2.3,1) circle (3pt);
\filldraw [color=blue] (0.5,1) circle (3pt);
\filldraw [color=blue] (0.5,2) circle (3pt);
\filldraw [color=blue] (0.5,3) circle (3pt);
\filldraw [color=blue] (0.5,4) circle (3pt);
\filldraw [color=orange] (2.7,2.5) circle (3pt);
\filldraw [color=green] (1.7,3) circle (3pt);
\filldraw [color=green] (1.7,4) circle (3pt);
\filldraw [color=green] (1.7,5) circle (3pt);
\filldraw [color=green] (1.7,6) circle (3pt);
\filldraw [color=orange] (3.5,3) circle (3pt);
\filldraw [color=red] (3,4) circle (3pt);
\filldraw [color=red] (3,5) circle (3pt);
\filldraw [color=red] (3,6) circle (3pt);
\filldraw [color=red] (3,7) circle (3pt);
\filldraw [color=blue] (4,0.5) circle (3pt);
\filldraw [color=blue] (4,2) circle (3pt);
\filldraw [color=orange] (4.5,3) circle (3pt);
\filldraw [color=red] (5,4) circle (3pt);
\filldraw [color=red] (5,5) circle (3pt);
\filldraw [color=red] (5,6) circle (3pt);
\filldraw [color=red] (5,7) circle (3pt);
\filldraw [color=blue] (5.7,1) circle (3pt);
\filldraw [color=orange] (5.5,2.5) circle (3pt);
\filldraw [color=green] (6.3,3) circle (3pt);
\filldraw [color=green] (6.3,4) circle (3pt);
\filldraw [color=green] (6.3,5) circle (3pt);
\filldraw [color=green] (6.3,6) circle (3pt);
\filldraw [color=orange] (7.5,1) circle (3pt);
\filldraw [color=orange] (7.5,2) circle (3pt);
\filldraw [color=orange] (7.5,3) circle (3pt);
\filldraw [color=orange] (7.5,4) circle (3pt);
\draw (2,0.6) node {\text{7}};
\draw (0.2,1) node {\text{4}};
\draw (0.2,2) node {\text{5}};
\draw (0.2,3) node {\text{8}};
\draw (0.2,4) node {\text{9}};
\draw (2.5,2.2) node {\text{21}};
\draw (1.6,3.3) node {\text{0}};
\draw (1.3,4) node {\text{1}};
\draw (1.3,5) node {\text{12}};
\draw (1.3,6) node {\text{13}};
\draw (3.2,2.7) node {\text{24}};
\draw (2.6,4) node {\text{18}};
\draw (2.6,5) node {\text{19}};
\draw (2.6,6) node {\text{30}};
\draw (2.6,7) node {\text{31}};
\draw (3.6,0.3) node {\text{6}};
\draw (3.6,1.7) node {\text{11}};
\draw (4.8,2.7) node {\text{25}};
\draw (5.3,4) node {\text{16}};
\draw (5.3,5) node {\text{17}};
\draw (5.3,6) node {\text{28}};
\draw (5.3,7) node {\text{29}};
\draw (6,0.7) node {\text{10}};
\draw (5.5,2.2) node {\text{20}};
\draw (6.3,3.3) node {\text{2}};
\draw (6.6,4) node {\text{3}};
\draw (6.6,5) node {\text{14}};
\draw (6.6,6) node {\text{15}};
\draw (7.8,1) node {\text{22}};
\draw (7.8,2) node {\text{23}};
\draw (7.8,3) node {\text{26}};
\draw (7.8,4) node {\text{27}};
\draw[->] (0.7,1) -- (2,1);
\draw[->] (0.7,1.9) -- (2.1,1.1);
\draw[->] (0.7,2.9) -- (2.2,1.2);
\draw[->] (0.7,3.9) -- (2.3,1.3);
\draw[->] (2.5,1) -- (3.8,0.5);
\draw[->] (1.8,2.8) -- (2.3,2.5);
\draw[->] (1.8,3.8) -- (2.4,2.6);
\draw[->] (1.8,4.8) -- (2.5,2.7);
\draw[->] (1.8,5.8) -- (2.6,2.8);
\draw[->] (2.9,2.4) -- (3.8,2);
\draw[->] (4,1.8) -- (4,0.7);
\draw[->] (5,3.8) -- (4.7,3);
\draw[->] (5,4.8) -- (4.6,3.1);
\draw[->] (4.9,5.8) -- (4.5,3.2);
\draw[->] (4.9,6.8) -- (4.4,3.3);
\draw[->] (4.4,2.8) -- (4.1,2.2);
\draw[->] (3,3.8) -- (3.3,3);
\draw[->] (3,4.8) -- (3.4,3.1);
\draw[->] (3,5.8) -- (3.5,3.2);
\draw[->] (3,6.8) -- (3.6,3.3);
\draw[->] (3.6,2.8) -- (3.9,2.2);
\draw[->] (6.2,2.9) -- (5.8,2.5);
\draw[->] (6.2,3.9) -- (5.7,2.6);
\draw[->] (6.2,4.8) -- (5.6,2.7);
\draw[->] (6.2,5.8) -- (5.5,2.8);
\draw[->] (5.3,2.4) -- (4.2,2);
\draw[->] (7.3,1) -- (6.1,1);
\draw[->] (7.3,1.9) -- (6,1.1);
\draw[->] (7.3,2.9) -- (5.9,1.2);
\draw[->] (7.3,3.9) -- (5.8,1.3);
\draw[->] (5.5,1) -- (4.3,0.5);
\draw[->] (3.9,0.3) .. controls (3.5,-0.15) and (4.5,-0.15) .. (4.1,0.3);
\filldraw [color=red] (10,5) circle (3pt);
\filldraw [color=green] (12,5) circle (3pt);
\filldraw [color=orange] (11,4) circle (3pt);
\filldraw [color=blue] (11,3) circle (3pt);
\draw (9.3,5) node {\text{(1,1)}};
\draw (12.7,5) node {\text{(1,-1)}};
\draw (11.7,4) node {\text{(-1,1)}};
\draw (10.2,3) node {\text{(-1,-1)}};
\draw[->] (10,4.8) -- (10.8,4.2);
\draw[->] (12,4.8) -- (11.2,4.2);
\draw[->] (11,3.8) -- (11,3.2);
\draw[->] (10.9,2.8) .. controls (10.5,2.35) and (11.5,2.35) .. (11.2,2.8);
\node at (7,7) {$\T_F$};
\draw[->] (7.5,7) -- (9.5,7);
\node at (8.5,7.5) {$\th$};
\node at (10,7) {$\T_\F$};
\end{tikzpicture}
\caption{
Above is a Boolean network with the dominant set $U=\{1\}$, and the corresponding induced logic network.
Below is the transition diagram of the Boolean network with states in $\tB^{\{1,\ldots,5\}}$ ordered lexicographically, and the transition diagram of the induced logic network. The eventual conjugacy $\th:\tB^{\{1,\ldots,5\}}\to\tB_2$ is codified by the correspondence of colors.}~\label{fig:Penacho}
\end{figure}
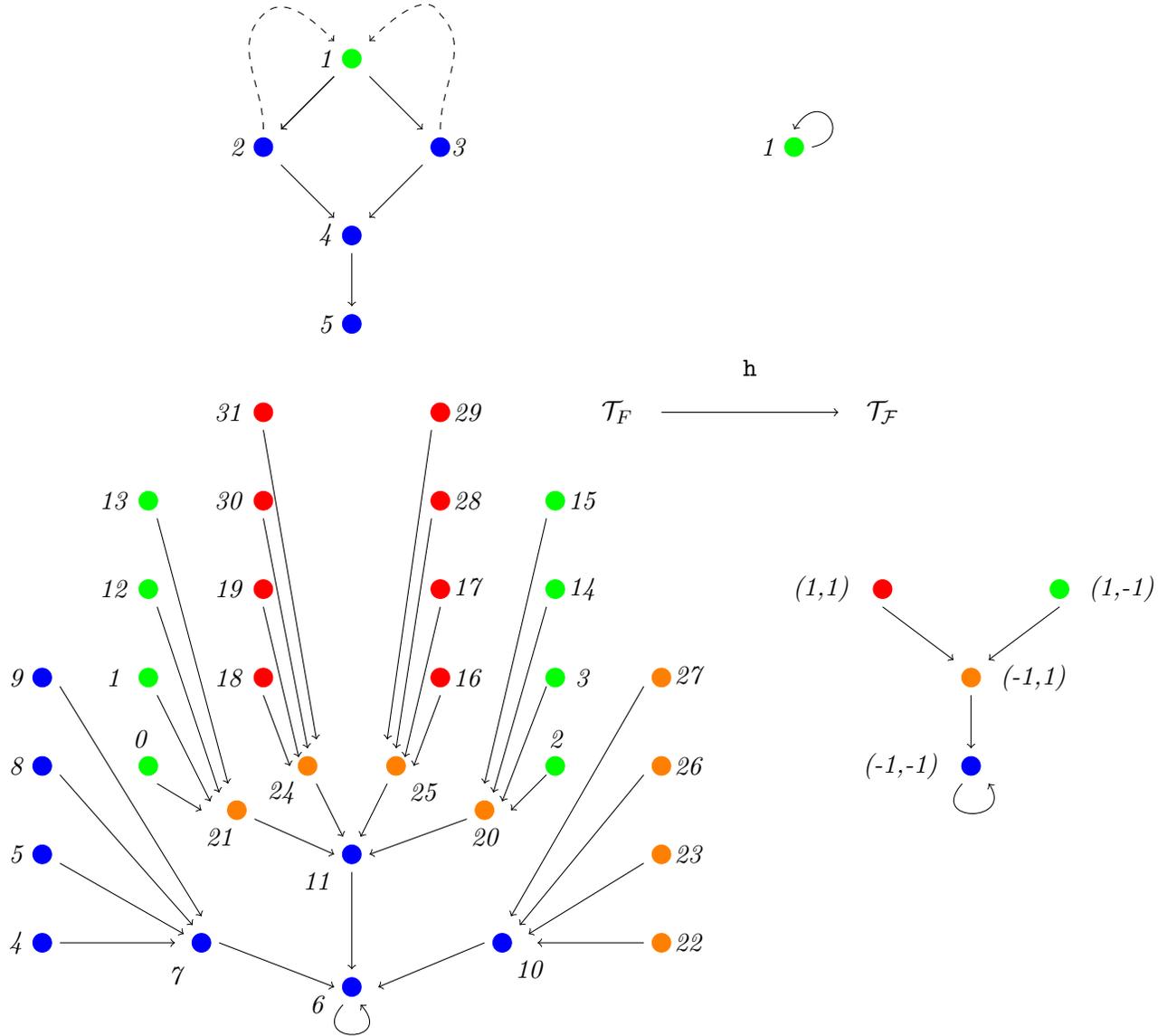
\end{center}

\end{example}

\ms As seen in the previous two examples, while the induced automata network preserves the structure of the attractors, it generally affects transients and the structure of the basins of attraction. Theorem~\ref{theo:Eventual-conjugacy} allows establishing bounds on the number of periodic attractors of each possible period, the depth of the transients, and the cardinality of the basins of attraction. Before stating these results, we will need some notation.

\ms Consider the Boolean network $(\tB^V,F)$ with underlying graph $(V,A)$. Let $U\subset V$ be a dominant set with depth $d=d(U)$, and recurrence length $\ell=\ell(U)$. Let $\left(\tB_\ell^U,\F\right)$ be the induced automata network, and $\th:\tB^V\to \tB_\ell^U$ the eventual conjugacy between both systems. Recall that ${\rm Per}_P(F)$ denotes the set of all $F$-periodic points of the period $P$ and ${\rm Per}(F)=\bigcup_{P\in\IN}{\rm Per}(F)$. For each $x\in (\tB^V, F)$, with $\tau_F(x)$ we denote the transient time of the orbit starting at $x$, that is, 
\[\tau_F(x):=\min\{t\in \INo:\, F^t(x)\in {\rm Per}(F)\}.\] 
As mentioned above, the landscape of the attractors comprises a finite number of basins attached to the periodic attractors, each corresponding to a connected component of the transition diagram $\T_F$. For each $x\in \tB^V$, denote by $C_F(x)$  the basin of attraction containing $x$. Finally, let us denote by $N_\A$ the number of attractors and by $N_\A(P)$ the number of those attractors of minimal period $P$. The same notions and corresponding notations apply to the induced automata network $(\tB_\ell^U,\F)$. We have the following results. 

\ms
\begin{corollary}[Attractors' landscape]~\label{cor:AttractorsLandscape}
Let $(\tB^V,F)$ be a Boolean network with underlying graph $(V, A)$. Let us suppose that $(V, A)$ admits a dominant set $U\subset V$ with depth $d=d(U)$ and recurrence length $\ell=\ell(U)$, and let $\left(\tB_\ell^U,\F\right)$ be the induced automata network, eventually conjugated to the original Boolean network through $\th:\tB^V\to \tB_\ell^U$. Then we have the following:
\begin{itemize}
\item[a)] $|{\rm Per}_P(F)|\leq 2^{p\times |U|}$ for each $P\in\IN$.
\item[b)] The maximal prime period, $P_{\max}:=\max\{P\in \IN: {\rm Per}_P(F)\setminus\bigcup_{P' < P}{\rm Per}_{P'}(F)\neq\emptyset\}$ is bounded by $2^{\ell\times |U|}$.
\item[c)] $N_\A\leq \sum_{P=1}^{P_{\max}}\frac{1}{P}\,2^{P\times |U|}$.
\item[d)] $|\tau_{\F}(\th(x))|\leq |\tau_F(x)|\leq |\tau_{\F}(\th(x))|+d$ for each $x\in \tB^V$.
\item[e)] $|C_{\F}(\th(x))|\leq |C_F(x)|\leq |C_{\F}(\th(x))|\times 2^{|V|-|U|}$ for each $x\in \tB^V$. 
\end{itemize}
\end{corollary}

\begin{proof}\

\begin{itemize}
\item[a)] Theorem~\ref{theo:Eventual-conjugacy} ensures that $|{\rm Per}_p(F)|=|{\rm Per}_p(\F)|$. If $y\in {\rm Per}_P(\F)$ and $p\leq \ell-1$, then 
\begin{align*}
y:=\left(y^0y^1\cdots y^{\ell-1}\right)&=\F^P(y)=\left(\Phi(\F^{p-1}(y))\cdots\Phi(y)\,y^0\cdots y^{\ell-P-1}\right),
\end{align*}
and therefore $y^t=y^{t+p}$ for $0\leq t\leq \ell-P-1$. Hence $y\in {\rm Per}_P(\F)$ is determined by its first $P$ coordinates $y^0,y^1,\ldots,y^{P-1}\in\tB^U$, which gives $|{\rm Per}_P(\F)|\leq 2^{P\times |U|}$. If on the contrary $P\geq \ell$, then $|{\rm Per}_P(\F)|\leq |\tB_\ell^U|= 2^{\ell\times|U|}\leq 2^{P\times|U|}$, and the claim follows.
\item[b)] It follows from the fact that 
\[x\in  {\rm Per}_P(F)\setminus\bigcup_{P' < P}{\rm Per}_{P'}(F)\text{ if and only if }
  h(x)\in {\rm Per}_P(\F)\setminus\bigcup_{P' < P}{\rm Per}_{P'}(\F),\] 
and in this case $|\{\F^t(h(x)):\, 0\leq t < P\}|=P\leq |\tB_\ell^U|= 2^{\ell\times|U|}$.  
\item[c)] It directly follows from a), b), and the fact that a periodic orbit of minimal period $p$ contains exactly $P$ points, hence
\[N_\A=\sum_{p=1}^{P_{\max}} N_P(F)=\sum_{P=1}^{2^\ell}N_P(\F)
=\sum_{P=1}^{P_{\max}} \frac{\left|{\rm Per}_P(\F)\setminus\bigcup_{P' < P}{\rm Per}_{P'}(\F)\right|}{P}
\leq \sum_{P=1}^{P_{\max}}\frac{2^{\ell\times|U|} }{P}.
\]
\item[d)] The inequality $|\tau_{\F}(\th(x))|\leq |\tau_F(x)|$ directly follows from the fact that $\th$ is a semi-conjugacy. 
Now, $\tau:=\tau_F(x)$ is by definition the smallest integer such that $F^\tau(x)=F^{\tau+P}(x)$, with $P\in \IN$ the eventual period of $x$. 
Similarly, $\F^{\tau'}(\th(x))=\F^{\tau'+P}(\th(x))$ for the first time when $\tau':=\tau_\F(\th(x))$. Notice that we have the same eventual period $P$, being $\th$ is an eventual conjugacy. 
Finally, since $\F^{\tau'}(\th(x))=\F^{\tau'+P}(\th(x))$, then $\F^{\tau'+s}(\th(x))=\F^{\tau'+P+s}(\th(x))$, i.e., $\th(F^{\tau'+s}(x))=\th(F^{\tau'+P+s}(x))$ for all $s\geq 0$. Therefore, $F^{\tau'+s}(x)_U=F^{\tau'+P+s}(x)_U$ for all $s\geq 0$, and invoking Theorem~\ref{theo:Dominance&Dynamics} we obtain $F^{\tau'+s}(x)_U=F^{\tau'+P+s}(x)_U$ for all $s\geq d$, with $d=d(U)$ the depth of the dominant set. We have in particular that $F^{\tau'+d}(x)_U=F^{\tau'+d+P}(x)_U$, and then necessarily $\tau\leq \tau'+d$.
\item[e)] Inequality $|C_{\F}(\th(x))|\leq |C_F(x)|$ follows directly from the fact that $\th$ is a semi-conjugacy and therefore $C_\F(\th(x))=\th(C_F(x))$. On the other hand, $ |C_F(x)|\leq |C_{\F}(\th(x))|\times 2^{|V|-|U|}$ is a consequence of the definition of $\th$, since $\th(x)=\th(x')$ implies that $x_U=x'_U$, and therefore $|\th^{-1}(y)|\leq 2^{|V|-|U|}$ for each $y\in\th(\tB^V)$.

\end{itemize}
\end{proof}

\ms
\paragraph{\bf Remarks to Corollary~\ref{cor:AttractorsLandscape}}\

\begin{enumerate}
\item Claim a) is related to Theorem 9 in~\cite{Aracena2008}, with the difference that his theorem gives the highest bound for the number of fixed points, which requires some knowledge of the local functions. In contrast, our claim depends only on the graph's topology and gives information about the number of periodic orbits of any period. 

\item It is not difficult to construct networks that trivially achieve the bounds established in Claim a). For this consider the cyclic graph on $p$ vertices: $0\mapsto 1\mapsto\cdots\mapsto (P-1)\mapsto 0$, with local functions $\phi_v(\zeta)=\zeta$ for each $\zeta\in\tB$. In this case $U=\{0\}$ is a dominant set and ${\rm Per}_P(F)=2^{p\times |U|}$. We obtain a system with $|U|$ arbitrary, considering a disjoint union of cycles as the one just described.

\item The maximal period $P=2^{\ell\times|U|}$ can be achieved. For the case $|U|=1$, consider the graph 
 $0\mapsto 1\mapsto\cdots\mapsto(\ell-1)$ and $v\mapsto 0$ for each $0\leq v < \ell$, with local functions $\phi_v(\zeta)=\zeta$ for $0< v < \ell$ and $\phi_0(x_0,x_1,\ldots,x_{\ell-1})$ such that the transitions 
\[(x_{\ell-1}x_{\ell-2}\cdots x_0)\rightarrow 
\left(x_{\ell-2}x_{\ell-3}\cdots x_0\phi_0(x_0,x_1\ldots\,x_{\ell-1})\right),\]
define an Eulerian path in the $\ell$-dimensional de-Bruijn graph in two symbols (see~\cite{deBruijn1946} for details). In this case, $U=\{0\}$ has recurrence length $\ell$, and the dynamics of the induced automata network consists of a single periodic orbit of period $2^{\ell\times|U|}$, visiting once each configuration in $\tB_\ell$. A similar construction can be devised for $U$ of any size.

\item The bounds in Claims d) and e) are achieved, i.e., there are networks for which $\tau_F(x)=\tau_\F(\th(x))$ and $|C_F(x)|=|C_\F(h(x))|\times 2^{|V|-|U|}$ for some $x\in\tB^V$. For $|U|=1$ and arbitrary $\ell$, the simplest example is obtained by considering the graph $\rotatebox[origin=c]{-90}{$\circlearrowright$}\,\, 0\mapsto 1\mapsto 2\mapsto\cdots\mapsto\,(|V|-1)$, and local functions $\phi_v(\zeta)=\zeta$ for each $0\leq v<|V|$.  The induced network has a single node and a local function $\Phi(y)=y$. In this case, $F$ has only two fixed-point attractors of the same size, hence $|C_F(x)|=2^{|V|-1}$ for each $x\in\tB^V$. On the other hand, $\F:\tB\to\tB$ has two fixed points, and therefore $|C_\F(y)|=1$ for each $y\in\tB$, and we have $|C_F(x)|=|C_\F(\th(x))|\times 2^{|V|-1}=|C_\F(\th(x)|\times 2^{|V|-|U|}$ for each $x\in\tB^V$ in this case.  Now, for $x=(-1)1^{|V|-1}$ we have $\tau_F(x)={|V|-1}$ while $\tau_\F(\th(x))=0$, therefore, $\tau_F(x)=\tau_\F(\th(x))+|V|-1=\F(\th(x))+d$.
\end{enumerate}

\ms
\section{Clover networks with signed majority rule}\label{sec:CloverNetworks}

\ms In this section we exemplify the reduction of the dynamics over a class of Boolean networks whose underlying graph $(V,A)$ have the following structure. The underlying graph $(V, A)$ has a distinguished vertex $v_0\in V$ connected via a directed path to each other vertex $v\in V$, and such that $|I(v)|=1$ for each $v\in V\setminus\{v_0\}$. Hence, the subgraph $(V, A\setminus I(v_0)\times \{v_0\})$ is a directed tree, rooted at $v_0$, and $(V, A)$ is completed by connecting each leave of this tree to $v_0$. Because of their shape, we will refer to graphs of this kind as clover networks.

\ms By construction, any clover network has $U=\{v_0\}$ as a dominant set, inducing an automata network on the loop $v_0\mapsto v_0$, and dynamics $\F:\tB_\ell\to\tB_\ell$, with $\ell-1=\max\{|C_k|:\, C_k \text{ is a cycle in } (V, A)\}$. For this class of Boolean networks, the induced automata network constitutes an effective reduction of the dynamics since $\ell\leq |V|$.  

\ms Concerning the local rules, we consider the following. We fix a map $\sigma:A\to \tB$, assigning an interaction sign to each $(u,v)\in A$. Then, for each $v\in V$, define $\phi_v:\tB^{I(v)}\to\tB$ by
\begin{equation}~\label{eq:Majority}
\phi_v\left(x_{I(v)}\right)={\rm Sign}\left(\left.\sum_{u\in I(v)}\sigma(u,v)x_u\right|\,x_v\right),
\end{equation}
where ${\rm Sign}:\IZ\times\tB\to\tB$ is such that
\[{\rm Sign}(y|\,x):=
\left\{\begin{array}{ll} {\rm Sign}(y) & \text{if }  y \neq 0,\\
                                    x  & \text{otherwise}.\end{array}\right.
\]         
A local rule of this kind is a signed majority rule, with signs given by the arrow map $\sigma$. In this case, to each directed cycle $C=v_0\mapsto v_1\mapsto\cdots\mapsto v_{m}\mapsto v_0$ in $(V,A)$, we associate the sign 
\begin{equation}\label{eq:CycleSign}
\sigma(C)=\sigma(v_m,v_0)\times \prod_{i=0}^{m-1}\sigma(v_i,v_{i+1}).
\end{equation}

\ms For clover networks with singed majority rule, we can explicitly compute the induced automata network. Indeed, we have the following.

\ms
\begin{proposition}[Induced from clover]~\label{prop:InducedClover} Let $(V,A)$ be a clover network and $\sigma:A\to \tB$ the interaction signs. Let $\S:=(\sum_{|C|=t}\sigma(C): 1\leq t\leq \ell)$ be the vector formed by the cycles' sign sums, grouped according to their length. Then, the corresponding Boolean network with signed majority rule induces an automata network such that $\F\left(y^0y^1\cdots y^{\ell-1}\right)=\left(\Phi(y)y^0y^1\cdots y^{\ell-1}\right)$, with $\Phi:\tB_\ell\to\tB$ such that
\[
\Phi\left(y^0y^1\ldots y^{\ell-1}\right)={\rm Sign}\left(\left.\sum_{t=0}^{\ell-1}\S_t\,y^{t-1}\right|\,y^0\right).
\]
\end{proposition}

\begin{proof}
By definition, 
\[\Phi\left(y^0y^1\cdots y^{\ell-1}\right)=\phi_{v_0}\left(x^0\right)={\rm Sign}\left(\left.\sum_{u\in I(v_0)}\sigma(u,v_0)x^0_u, x_{v_0}\right|\, x_{v_0}\right),\] 
with $x^0\in \tB^{I(v_0)}$ defined as follows.
\begin{itemize}
\item[] If $v_0\in I(v_0)$, then $x^0_{v_0}=y^0$. In this case $C_{v_0}:=v_0\mapsto v_0$ is a cycle with $\sigma(C)=\sigma(v_0,v_0)$,
\item[] For $u\in I(v_0)\setminus \{v_0\}$, $x^0_{u}=\phi_{u}\left(x^1_v\right)={\rm Sign}\left(\sigma(v,u)x_v^1\right)=\sigma(v,u)x_v^1$,  where $I(u)=\{v\}$. 
\end{itemize}
For each $t\geq 1$, $x^t\in \tB^{I^{t+1}(v_0)}$ is defined recursively:
\begin{itemize}
\item[] If $v_0\in I^t(v_0)$, then $x^t_{v_0}=y^t$. In this case, $v_0$ appears in a cycle of length $t+1$.
\item[] For $u\in I^t(v_0)\setminus\{v_0\}$, $x^t_{u}=\phi_{u}\left(x^{t+1}_v\right)={\rm Sign}\left(\sigma(v,u)x_v^{t+1}\right)=\sigma(v,u)x^{t+1}_v$,  where $I(u)=\{v\}$.
\end{itemize}
The recursion stops at $t=\ell-1$, with $\ell=\max\{|C|:\, C \text{ is a cycle}\}$.
The clover structure ensures that each $u\in I(v_0)$ appears only in one cycle $C_u=v_0\mapsto \cdots\mapsto u\mapsto v_0$, and for that cycle, the recursion above gives $x_u^0=\sigma(C_u)y^{|C_u|-1}$, from which the proposition follows. 
\end{proof}

\ms
\paragraph{\bf Remarks to Proposition~\ref{prop:InducedClover}}\

\begin{enumerate}
\item The dynamics $F: \tB^V\to \tB^V$, defined by Equation~\eqref{eq:Majority}, commutes with the symmetry $x\mapsto -x$ in $\tB^V$. This derives from the fact that ${\rm Sign}(-x|\,-y)=-{\rm Sign}(x,y)$ for each $x\in \IZ,\, y\in \tB$. The induced automata, $\F:\tB_\ell\to\tB_\ell$ inherits this symmetry.
\item The symmetry $x\mapsto -x$ reflects in a symmetry in the transition diagram for both, the original and the induced automatq.
\item A sufficient condition for two clover networks to be eventually conjugated is that they give place to the same vector $\S:=(\sum_{|C|=t}\sigma(C): 1\leq t\leq \ell)$, which groups the sum of cycle signs of the same length.
\item Suppose that $\S$ is such that $|\S_\tau| > \sum_{t\neq \tau}^{\ell}|\S_\tau|$ for some $1\leq \tau \leq \ell$. In this case the induced automata network is such that $y\mapsto\left({\rm Sign}(\S_\tau\, y^{\tau-1})y^0\cdots y^{\ell-2}\right)$ and we have two possibilities: 
\begin{itemize}
\item[a)] For $\S_\tau > 0$, each initial condition $y\tB_\ell$ converges to a periodic attractor of period $\tau$, completely determined by the prefix $y^0y^1\cdots y^{\tau-1}$ of $y\in \tB_\ell$.
\item[b)] If $\S_\tau <0$, then each initial condition $y\tB_\ell$ converges to a periodic attractor of period $2\,\tau$, determined as well, by the prefix $y^0y^1\cdots y^{\tau-1}$. 
\end{itemize}
More generally, the attractors' landscape can be easily determined when a few coordinates of $\S$ dominate over the rest.
\end{enumerate}

\ms 
\begin{example}
Consider the Boolean network $\left(\tB^{\{1,\ldots,5\}},F\right)$ depicted in Figure~\ref{fig:Trebol}, with signed majority local rules as defined by Equation~\eqref{eq:Majority}, using the interaction signs indicated in the same figure. For this clover network,  the dominant set is $U=\{1\}$, with depth $d=1$, and recurrence length $\ell=2$. Hence, the induced automata network is a transformation $\F:\tB_2\to\tB_2$, defined by the local function
\begin{equation}~\label{eq:InducedTrebol}
\Phi\left(y^0\,y^1\right)={\rm Sign}\left(-2\,y^1|y^0\right)=-{\rm Sign}\left(y^1\right).
\end{equation}
In the same figure, we show the transition diagrams of both the original Boolean network and the induced automata network. The eventual conjugacy $\th:\tB^{\{1,\ldots,5\}}\to\tB_2$ is codified by the correspondence of node colors in both diagrams.

\begin{center}
\begin{figure}
\begin{tikzpicture}[scale=1.3]
\filldraw [color=green] (4,10) circle (3pt);
\filldraw [color=blue] (3,10) circle (3pt);
\filldraw [color=blue] (5,10) circle (3pt);
\filldraw [color=blue] (4,11) circle (3pt);
\filldraw [color=blue] (4,9) circle (3pt);
\draw (4,9.7) node {\text{1}};
\draw (5.2,10) node {\text{2}};
\draw (3.7,11) node {\text{3}};
\draw (2.7,10) node {\text{4}};
\draw (3.7,9) node {\text{5}};
\filldraw [color=green] (9,10) circle (3pt);
\draw (8.7,10) node {\text{1}};
\draw[->] (4.2,9.9)  to [out=-45,in=-135]  (4.8,9.9);
\draw[->] (3.2,9.9)  to [out=-45,in=-135]  (3.8,9.9);
\draw[->] (4.8,10.1)  to [out=115,in=45]  (4.2,10.1);
\draw[->] (3.8,10.1)  to [out=155,in=45]  (3.2,10.1);
\draw[->] (4.1,10.2)  to [out=45,in=-45]  (4.1,10.8);
\draw[->] (3.9,10.8)  to [out=-155,in=155]  (3.9,10.2);
\draw[->] (4.1,9.2)  to [out=45,in=-45]  (4.1,9.8);
\draw[->] (3.9,9.8)  to [out=-155,in=155]  (3.9,9.2);
\draw[->] (9.2,10) .. controls (9.7,10.1) and (9.3,10.7) .. (9,10.2);

\draw (4.7,9.6) node {\text{--}};
\draw (4.7,10.4) node {\text{+}};
\draw (4.3,10.7) node {\text{--}};
\draw (3.6,10.7) node {\text{+}};
\draw (3.3,10.4) node {\text{+}};
\draw (3.3,9.6) node {\text{+}};
\draw (4.3,9.3) node {\text{--}};
\draw (3.6,9.3) node {\text{+}};


\filldraw [color=red] (3,5) circle (3pt);
\filldraw [color=green] (5,5) circle (3pt);
\filldraw [color=orange] (3,4) circle (3pt);
\filldraw [color=blue] (5,4) circle (3pt);
\draw (3.2,4.7) node {\text{3}};
\draw (4.7,4.7) node {\text{12}};
\draw (3.3,4.3) node {\text{19}};
\draw (4.7,4.3) node {\text{28}};
\draw[->] (3.2,5) -- (4.8,5);
\draw[->] (4.8,4) -- (3.2,4);
\draw[->] (3,4.2) -- (3,4.8);
\draw[->] (5,4.8) -- (5,4.2);

\filldraw [color=orange] (2,5) circle (3pt);
\filldraw [color=orange] (2,5.5) circle (3pt);
\filldraw [color=orange] (2,6) circle (3pt);
\filldraw [color=orange] (2.5,6) circle (3pt);
\filldraw [color=orange] (3,6) circle (3pt);
\draw (1.6,5) node {\text{16}};
\draw (1.6,5.5) node {\text{17}};
\draw (1.6,6) node {\text{18}};
\draw (2.5,6.3) node {\text{21}};
\draw (3,6.3) node {\text{25}};
\draw[->] (2.2,5) -- (2.8,5);
\draw[->] (2.2,5.4) -- (2.8,5.1);
\draw[->] (2.2,5.9) -- (2.8,5.2);
\draw[->] (2.5,5.8) -- (2.9,5.2);
\draw[->] (3,5.8) -- (3,5.2);

\filldraw [color=green] (6,4) circle (3pt);
\filldraw [color=green] (6,3.5) circle (3pt);
\filldraw [color=green] (6,3) circle (3pt);
\filldraw [color=green] (5.5,3) circle (3pt);
\filldraw [color=green] (5,3) circle (3pt);
\draw (6.3,4) node {\text{6}};
\draw (6.3,3.5) node {\text{10}};
\draw (6.3,3) node {\text{13}};
\draw (5.5,2.7) node {\text{14}};
\draw (5,2.7) node {\text{15}};
\draw[->] (5.8,4) -- (5.2,4);
\draw[->] (5.8,3.5) -- (5.2,3.9);
\draw[->] (5.8,3.2) -- (5.2,3.8);
\draw[->] (5.5,3.2) -- (5.1,3.8);
\draw[->] (5,3.2) -- (5,3.8);

\filldraw [color=red] (3.5,6) circle (3pt);
\filldraw [color=red] (4,6) circle (3pt);
\filldraw [color=red] (4.5,6) circle (3pt);
\filldraw [color=red] (5,6) circle (3pt);
\filldraw [color=red] (5.5,6) circle (3pt);
\filldraw [color=red] (6,6) circle (3pt);
\filldraw [color=red] (6,5.5) circle (3pt);
\filldraw [color=red] (6,5) circle (3pt);
\filldraw [color=red] (6,4.5) circle (3pt);
\draw (3.5,6.3) node {\text{0}};
\draw (4,6.3) node {\text{1}};
\draw (4.5,6.3) node {\text{2}};
\draw (5,6.3) node {\text{4}};
\draw (5.5,6.3) node {\text{5}};
\draw (6,6.3) node {\text{7}};
\draw (6.3,5.5) node {\text{8}};
\draw (6.3,5) node {\text{9}};
\draw (6.3,4.5) node {\text{11}};
\draw[->] (3.5,5.8) -- (4.8,5.1);
\draw[->] (4,5.8) -- (4.8,5.2);
\draw[->] (4.5,5.8) -- (4.9,5.2);
\draw[->] (5,5.8) -- (5,5.2);
\draw[->] (5.5,5.8) -- (5.1,5.2);
\draw[->] (5.8,5.8) -- (5.2,5.2);
\draw[->] (5.8,5.5) -- (5.2,5.1);
\draw[->] (5.8,5) -- (5.2,5);
\draw[->] (5.8,4.5) -- (5.2,4.9);

\filldraw [color=blue] (4.5,3) circle (3pt);
\filldraw [color=blue] (4,3) circle (3pt);
\filldraw [color=blue] (3.5,3) circle (3pt);
\filldraw [color=blue] (3,3) circle (3pt);
\filldraw [color=blue] (2.5,3) circle (3pt);
\filldraw [color=blue] (2,3) circle (3pt);
\filldraw [color=blue] (2,3.5) circle (3pt);
\filldraw [color=blue] (2,4) circle (3pt);
\filldraw [color=blue] (2,4.5) circle (3pt);
\draw (4.5,2.7) node {\text{20}};
\draw (4,2.7) node {\text{22}};
\draw (3.5,2.7) node {\text{23}};
\draw (3,2.7) node {\text{24}};
\draw (2.5,2.7) node {\text{26}};
\draw (2,2.7) node {\text{27}};
\draw (1.6,3.5) node {\text{29}};
\draw (1.6,4) node {\text{30}};
\draw (1.6,4.5) node {\text{31}};
\draw[->] (4.4,3.2) -- (3.2,3.9);
\draw[->] (4,3.2) -- (3.15,3.8);  
\draw[->] (3.5,3.2)  -- (3.05,3.8);
\draw[->] (3,3.2) -- (3,3.8);
\draw[->] (2.5,3.2) -- (2.95,3.8);
\draw[->] (2.1,3.2) -- (2.9,3.85);
\draw[->] (2.2,3.5) -- (2.8,3.9);
\draw[->] (2.2,4) -- (2.8,4);
\draw[->] (2.2,4.5) -- (2.8,4.1);


\filldraw [color=red] (10,5) circle (3pt);
\filldraw [color=green] (12,5) circle (3pt);
\filldraw [color=orange] (10,4) circle (3pt);
\filldraw [color=blue] (12,4) circle (3pt);
\draw (9.3,5) node {\text{(-1,-1)}};
\draw (12.7,5) node {\text{(1,-1)}};
\draw (9.3,4) node {\text{(-1,1)}};
\draw (12.7,4) node {\text{(1,1)}};
\draw[->] (10.2,5) -- (11.8,5);
\draw[->] (11.8,4) -- (10.2,4);
\draw[->] (10,4.2) -- (10,4.8);
\draw[->] (12,4.8) -- (12,4.2);
\node at (7,7) {$\T_F$};
\draw[->] (7.5,7) -- (9.5,7);
\node at (8.5,7.5) {$\th$};
\node at (10,7) {$\T_\F$};
\end{tikzpicture}
\caption{
Above is a clover network with signed majority rule defined by the sings depicted next to the arrows. The the corresponding induced automata network is a sigle loop at vertex 1.
Below is the transition diagram of the original network with states in $\tB^{\{1,\ldots,5\}}$ ordered lexicographically, and the transition diagram of the induced automata network. The eventual conjugacy $\th:\tB^{\{1,\ldots,5\}}\to\tB_2$ is codified by the correspondence of colors.}~\label{fig:Trebol}
\end{figure}
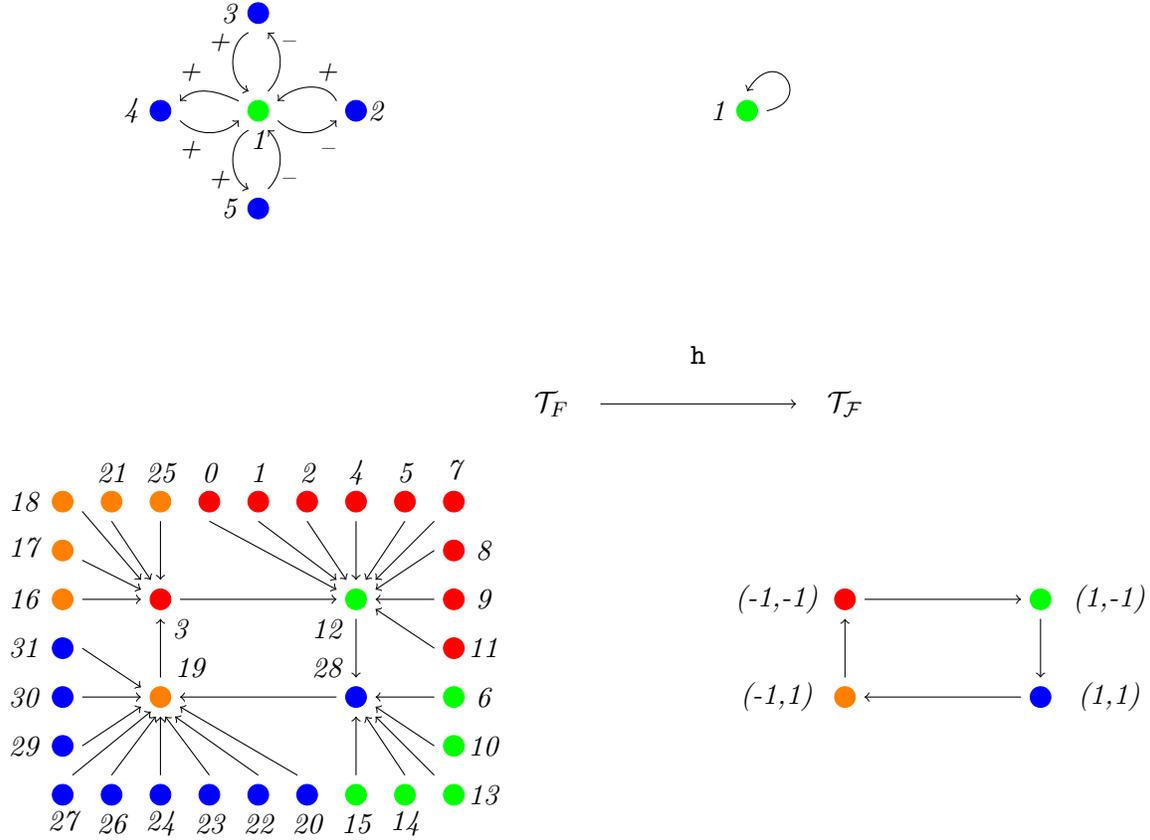
\end{center}    

\end{example}

\ms In what follows we analyze the landscape of attractors in an ensemble of clover Boolean networks with signed majority rule. We compare the mean structure of the attractor's landscape of the original Boolean network to that of induced automata network. The aim is to evaluate the influence of the sign distribution and the distribution of cycle lengths on key dynamical features such as the number of attractors, their periods, the size of their respective basins, and the transient behavior preceding convergence. At the same time we will study the effect of network reduction, on the transient behavior. In this way we can asses who loose are the theoretical bounds established in Corollary~\ref{cor:AttractorsLandscape}, concerning the transient lengths and the size of the basins of attraction.  

\ms
\paragraph{\bf An ensemble of clover networks with signed majority rule}
For $N\in\IN$, $p\in(0,1)$ and $q\in [0,1]$ we construct a clover network with signed majority rule as follows: 
\begin{itemize}
\item[a)] For each $2< n\leq N$, put $(1,n)\in A$ uniformly and independently with probability $p$. In this way we obtain the output set $O(1)=\{n_0:=2,n_1,n_2,\ldots,n_s\}\subset \{2,3,\ldots,N\}$, such that $(1,n)\in A$ for each $n\in O(1)$.  
\item[b)] Let $n_{s+1}:=N+1$ and for each $0\leq k\leq s$ and $n_k\leq n < n_{k+1}-1$, let $(n,n+1)\in A$.
\item[c)] For each $1\leq k\leq s+1$, $(n_k-1,1)\in A$.
\end{itemize}   
In this way we obtain a clover network with cycles $1\mapsto n_k\mapsto (n_k+1)\mapsto\cdots\mapsto (n_{k+1}-1)\mapsto 1$, for $1\leq k\leq s+1$. The parameter $p$ controls the number and the size of the of cycles of the clover network, and we will refer to it as the folding probability. Indeed, the number $N_c$ of cycles is distributed according to $\IP(N_c=k)=1+f(k,N-2,p)$, with $f(k,N-2,p)$ the binomial distribution with $N-2$ trials and parameter $p$. Hence, the expected number of cycles is $\IE(N_c)=1+p(N-2)$. The length $\ell_1$ of the cycle $C_1$ is such that  $\IP(\ell_1=\ell)=p^{\chi_{<N}(\ell)}\,(1-p)^{\ell-2}$, with  $\chi_{<N}$ is the characteristic function of the set $\{\ell<N\}$. This gives an expected length $\IE(\ell_1)=\left(-p^2 -(1-p)^N\right)/\left(p(1-p)\right)\rightarrow (1+p)/p$, as $N\to\infty$. This asymptotically coincides with the common length of the cycles in a clover network with $1+p(N-2)$ cycles of the same length. We can estimate the distribution of the maximal length $\max_{1\leq t\leq \ell}\ell_t$ by considering that our random structure is reasonably modeled by a clover network with $N_c=1+p(N-2)$ cycles whose lengths are i.i.d. random variables with exponential distribution $\IP(\ell_k=\ell)=p(1-p)^{\ell-2}$. For this simplified model, 
\[
\IP\left(\max_{1\leq k\leq N_c}\ell_k \geq \ell\right)=1-\left(1-(1-p)^{\ell-2}\right)^{N_c},
\]  
and therefore  
\begin{eqnarray}~\label{eq:ExpecMax}
\IE\left(\max_{1\leq k\leq N_c}\ell\right)
    &=& 1 + \sum_{\ell=2}^\infty\IP\left(\max_{1\leq k\leq N_c}\ell_k \geq \ell\right)
     = 1 + \sum_{m=0}^\infty\left(1-\left(1-(1-p)^{\ell-2}\right)^{N_c}\right)\\ \nonumber
    &=&1 + \sum_{m=0}^\infty(1-p)^m\left(\sum_{k=0}^{N_c-1}\left(1-(1-p)^{m}\right)^k\right) \\ \nonumber
    &\approx & 1+\sum_{k=0}^{N_c-1}\int_0^\infty (1-p)^m(1-(1-p)^m)^k\, dm 
    = 1 + \frac{\sum_{k=1}^{N_c}\left(1-p^k\right)\,k^{-1}}{-\log(1-p)} \\ \nonumber
    &\sim & \frac{\log(N_c)+\gamma-(1-p)^{-1}}{-\log(1-p)}
     \sim  \frac{\log(N)+\log(p)+\gamma-(1-p)^{-1}}{-\log(1-p)},
\end{eqnarray}
with $\gamma\approx 0.57721$ the Euler-Mascheroni constant. This estimation agrees with numerical estimations on our original ensemble, and its asymptotical validity can be rigorously established.  

\ms To the network so obtained we associate interaction signs by maps $\sigma:A\to\tB$ determined as follows. For $(i,j)\in A$, let $\sigma(i,j)=-1$ uniformly and independently with probability $q\in [0,1]$. The parameter $q$, which we call the inhibition probability, determines the the distribution of cycle signs. Since for the clover networks that we consider, cycles intersect only at vertex $1$, then the distribution of the cycle sings are independent. It can be easily verified that, for a cycle of $C$ of length $L$, we have
\begin{equation}~\label{eq:etaLq}
\IP(\sigma(C)=1)=\eta_L(q):=\sum_{m=0}^{\lfloor L/2\rfloor}\left(\begin{matrix}
N\\2m
\end{matrix}\right)q^{2m}(1-q)^{L-2m}=\frac{(1-2q)^L+1}{2}. 
\end{equation}
Notice that $\eta_L(q)\rightarrow 1/2$ as $L\to\infty$, for each $0<q<1$. 

\ms Given the previous estimations, we expect the networks in the ensemble to satisfy the following:
\begin{enumerate}
\item According to Equation~\eqref{eq:ExpecMax}, the expected recurrence length of a clover network in our ensemble grows as $O\left(\log(N)\right)$, hence, by Corollary~\ref{cor:AttractorsLandscape}-b), the maximal expected prime period is bounded by a power-law $N\mapsto \alpha_p\, N^{\beta_p}$ with $\beta_p=-\log(2)/\log(1-p)$. 
\item Since the distribution of cycles is approximately exponential, the majority of cycles would be of short length. In this case, the vector $\S:=(\sum_{|C|=t}\sigma(C): 1\leq t\leq \ell)$, which groups the sum of cycle signs of the same length, is expected to be dominated by its first coordinates. In this case, the maximal expected period of the induced automata network would be the number of those first dominating coordinates, and the transients as long as the number of complementary coordinates. The number of dominating coordinates is expected to be a decreasing function of the folding parameter.
\item According to Equation~\eqref{eq:etaLq}, the influence of the inhibition probability becomes less important as the size of the cycles grows. Hence, the large coordinates of $\S$ can be modeled by centered random variables with a variance tending to zero. 

\end{enumerate}

\ms
\paragraph{\bf Numerical exploration of the induction effect on the attractors' landscape}
For fixed $N,p,q$, we characterize the structure of the transition diagrams $\T_F$ and $\T_\F$ by considering the following indicators:
\begin{itemize}
\item[a)] The common number $N_\A$ of connected components of $\T_F$ and $\T_\F$, i.e., the common number of attractors. 
\item[b)] The reduction in size, $|\C_{\mathcal{F}}|/|\C_F|$, of each connected component after induction, i.e., the cardinality ratio of the corresponding basins of attraction.
\item[c)] The common length $P(\C)$ of the cyclic core of each component of $\T_F$ and $\T_\F$, i.e., the common period of the corresponding periodic attractors.
\item[d)] The reduction after induction of the mean depth $\Delta\tau=\langle \tau_F\rangle_{C_F}-\langle \tau_\F\rangle_{\C_\F}$ of each component, which is noting but the difference between the basins' mean transients before and after induction.  
\item[e)] The reduction after induction of the maximal depth $\Delta\tau_{\max}=\max_{x\in C}\tau_F(x)-\max_{y\in\th(C)}\tau_\F(y)$ of each component, which is noting but the difference between the basins' maximal transients before and after induction.  
\end{itemize}
We register as well, the recurrence length $\ell$, of the underlying clover network, which serves a key parameter for the theoretical upper bounds established in Corollary~\ref{cor:AttractorsLandscape}.

\ms We performed 500 simulations with $N=10$, for each $p,q \in \{0.3,0.6,0.9\}$, observing some interesting features regarding the attractors' landscape. Notice that the total number of admissible clover configurations is bounded above by $2\cdot 3^{N-2}=2\cdot 3^8$. In the case of normally distributed observables, this number of realizations ensures, in the worst-case scenario, a confidence level of $95\%$ with a margin of error below $5\%$ for proportion-type observables.

\begin{enumerate}
\item The number $N_\A$ of attractors was, in all of our simulations, much smaller than the theoretical upper bound $\sum_{P=1}^{P_{\max}}2^P/P$, which indicates that the number of attractors of a given period, $N_\A(P)$, is usually smaller than the theoretical upper bound $2^P/P$. Besides, $N_\A$ seems to decrease with the folding probability $p$, consistent with the fact that $N_\A$ increases with the recurrence length $\ell$. For small $p$, the distribution of $N_\A$ is highly skewed and exhibits a wide $[q_{0.05},q_{0.95}]$ interval, which captures the central 90\% of the realizations and therefore reflects substantial dispersion. For larger $p$ the interval contracts significantly, indicating a clear concentration phenomenon. 

\ms
\begin{center}
\begin{table}[h]
\begin{tabular}{|l|c|c|c|}
\hline
\backslashbox{$q$}{$p$} & 0.3 & 0.6 & 0.9\\
\hline
0.3 & 26.5\, [1,52] & 3.0\, [1,5] & 2.0\, [1,3]\\
\hline
0.6 & 22.5\, [1,44] & 3.5\, [1,6] & 2.0\, [1,3]\\
\hline
0.9 & 19.0\, [2,36] & 3.5\, [2,5] & 3.0\, [3,3]\\
\hline
\hline
$\ell$						& 5.87 $\pm$ 1.71	&4.06	$\pm$ 1.07 &2.66 $\pm$ 0.67\\
\hline
\end{tabular}
\caption{The median number of attractors, $N_\A$, together with the empirical $[q_{0.05},q_{0.95}]$ interval over 500 simulations, as a function of the folding and inhibition probabilities.  In the last line, the recurrence length $\ell$ is reported as mean $\pm$ standard deviation.  All simulations were performed with $N=10$ nodes.}
\end{table}
\end{center}

\item The distribution of periods along the attractors is heavily skewed toward short periods, with a clearly non-normal shape characterized by a decay with a thick tail, and although the theoretical upper bound grows exponentially with $\ell$, the observed periods were always much smaller. This is consistent with the fact that the distribution of cycles in a clover network generates an induced local function $\Phi:\tB_\ell\to\tB$ effectively depending on a few first coordinates. Indeed, in all of our simulations, the maximal observed period was upper bounded by $2\ell$. The distribution of periods appears to depend on both the folding probability $p$ and the inhibition probability $q$, but most notably on $p$. The mean period $\langle P\rangle:=\frac{1}{N_\A}\sum_\C P(\C)$, appears to be a decreasing function of $p$, and for each $p$ fixed, it would be a unimodal function of $q$ with maximum presumably around $q=1/2$. 

\ms
\begin{center}
\begin{table}[h]
\begin{tabular}{|l|c|c|c|}
\hline
\backslashbox{$q$}{$p$} & 0.3 & 0.6 & 0.9\\
\hline
0.3 & 7.3\, [1.0,13.6] & 4.5\, [1.0,8.0] & 2.5\, [1.0,4.0]\\
\hline
0.6 & 7.5\, [1.3,13.6] & 4.5\, [1.0,8.0] & 2.6\, [1.3,4.0]\\
\hline
0.9 & 6.1\, [1.3,11.0] & 3.6\, [1.3,6.0] & 1.3\, [1.3,1.3]\\
\hline
\end{tabular}
\caption{Median of the mean period $\langle P\rangle$ (computed per realization) together with the empirical $[q_{0.05},q_{0.95}]$ interval over 500 simulations, as a function of the folding and inhibition probabilities.}
\end{table}
\end{center}

\item The reduction in size of each basins of attraction, $|\C_{\mathcal{F}}|/|\C_F|$, is rather uniform and close to the quotient $\tB_\ell/\tB^N=2^{\ell-N}$. This reduction factor remains far from the theoretical bound $2^{1-N}$, which would correspond to a very biased distribution of reductions among the different attractors.

\ms
\begin{center}
\begin{table}[h]
\begin{tabular}{|l|c|c|c|}
\hline
\backslashbox{$q$}{$p$} & 0.3 & 0.6 & 0.9\\
\hline
0.3 & 0.507\, [0.015,1.000] & 0.035\, [0.007,0.062] & 0.009\, [0.003,0.015]\\
\hline
0.6 & 0.382\, [0.015,0.750] & 0.035\, [0.007,0.062] & 0.009\, [0.003,0.015]\\
\hline
0.9 & 0.257\, [0.015,0.500] & 0.035\, [0.007,0.062] & 0.009\, [0.003,0.015]\\
\hline
\hline
$2^{-(N-\ell)}$ & 0.057 $\pm$ 0.067 & 0.016 $\pm$ 0.012 & 0.006 $\pm$ 0.002\\
\hline
\end{tabular}
\caption{Median of the reduction factor $|\C_{\mathcal{F}}|/|\C_F|$ over 500 realizations, together with the empirical $[q_{0.05},q_{0.95}]$ interval. The last line shows the theoretical mean reduction $2^{-(N-\ell)}$, with $\ell$ reported as mean $\pm$ standard deviation. The comparatively large dispersion induced in $2^{-(N-\ell)}$ is a direct consequence of its exponential dependence on $\ell$, and does not indicate statistical instability of the sample.}
\end{table}
\end{center}

\item The reduction in mean and maximal transient length appears to be quite uniform as well. This is even more striking in the case of the maximal transient length, $\Delta\tau_{\max}$, which we found to be independent of attractor in more than 90\% of our simulations. Both $\Delta\tau$ and $\Delta\tau_{\max}$ were upper bounded by $\ell-1$ in all of our simulations, in accordance with our theoretical upper bound $\tau_F(x)-\tau_\F(\th(x))\leq d=\ell-1$. Both differences appear to be decreasing with $p$ and $q$, indicating that for $p$ and $q$ approaching 1, the transients are less reduced. The upper bound $\ell-1$ itself decreases with $p$, in agreement with the estimate for clover networks. 

\ms
\begin{center}
\begin{table}[h]
\begin{tabular}{|l|c|c|c|}
\hline
\backslashbox{$q$}{$p$} & 0.3 & 0.6 & 0.9\\
\hline
0.3 & 2.0\, [0.0,4.0] & 2.0\, [1.0,3.0] & 1.5\, [1.0,2.0]\\
\hline
0.6 & 2.1\, [0.3,4.0] & 2.0\, [1.0,3.0] & 1.5\, [1.0,2.0]\\
\hline
0.9 & 2.3\, [0.6,4.0] & 1.8\, [1.0,2.6] & 1.0\, [1.0,1.0]\\
\hline
\hline
$\ell$ & 5.87 $\pm$ 1.71 & 4.06 $\pm$ 1.07 & 2.66 $\pm$ 0.67\\
\hline
\end{tabular}
\caption{Median of the reduction in maximal transient length $\Delta\tau_{\max}$ over 500 realizations, together with the empirical $[q_{0.05},q_{0.95}]$ interval, as a function of the folding and inhibition probabilities. In the last line, the recurrence length $\ell$ is reported as mean $\pm$ standard deviation. All simulations were performed with $N=10$.}
\end{table}
\end{center}

\end{enumerate}

\section{Concluding Remarks}\label{sec:Conclusions}

\ms Given a set of dominant vertices in a Boolean network under a synchronous update scheme, we induce an automata network on those vertices, which is asymptotically equivalent to the original one. In this setting, when the recurrence length and the cardinality of the dominant set are small, the induced automata network constitutes a reduction of the original Boolean network, which preserves all the periodic attractors while reducing the corresponding basins. This reduction depends only on the topology of the original network and is valid only in the case of synchronous updating. 

\ms The network reduction, when possible, could be further improved by taking into account the particularities of the regulatory rules. This is observed in the case of clover networks with signed majority rule, where the induced local rule depends only on a few first coordinates. A refinement of this reduction, by taking into account the description of the local rules, will be considered in future works. Since dominant sets of vertices make sense for asynchronous updates, similar results are expected for the asynchronous update schemes. Nevertheless, the extension to that framework would require a different treatment of dependency propagation. Such an extension is of particular interest due to its potential impact on the dynamical analysis of bio-inspired networks and related discrete dynamical systems, and is therefore left for future work.

\ms The proposed reduction does not provide an \emph{a priori} estimate of its size, since the cardinality of dominant sets depends entirely on the network topology. In the examples presented here, the reduction can be extreme and the dynamics of a network with $N$ vertices is fully captured by a single dominant vertex, yielding an effective complexity scaling of order $1/N$. Typical reduction factors in biologically inspired networks are currently under investigation.

\ms We consider that the notion of eventual equivalence, which we introduced in this paper, deserves further study. It could be used to classify dissipative dynamical systems, among other possible applications. In~\cite{CoutinhoEatl2006}, the symbolic description of the asymptotic dynamics of piecewise contractions leads to associating a dissipative system with some low-complexity dynamical systems, as the irrational rotations. These symbolic descriptions of the asymptotic dynamics should be considered in the framework of an asymptotic equivalence.    

\ms Our results suggest that, given the interactions' structure, a maximal possible complexity of the dynamics is encoded in the induced network defined by the dominant vertices. An interesting program would consist of studying the landscape of attractors for induced networks of increasing structural complexity, which would be the simplest representatives of each class of asymptotically equivalent systems. 

\ms In future work we will focus on two main directions:  (i) study of the induced automata network from gene regulatory networks, where dominant vertices may capture biologically meaningful control modules and simplify the analysis of cell-fate decision processes; and (ii) extensive statistical studies of ensembles of random networks, including Erd\"os-R\'enyi and Barab\'asi-\'Albert topologies, to investigate how structural randomness and degree heterogeneity affect the landscape of attractors and how this correlates with the structure of dominant vertices for those ensembles.

\section*{Code Availability}

\ms The implementation used to generate clover networks and compute their full and reduced dynamics is publicly available at:
\url{https://github.com/arletteespana/clover_dynamics}

\section*{Acknowledgments}

\ms The authors express their gratitude to the Secretar\'{\i}a de Ciencia, Humanidades, Tecnolog\'{\i}a e Innovaci\'on (SECIHTI) for the institutional support and the promotion of research that made this work possible.

\bibliographystyle{plain}

\end{document}